\documentclass[12pt,a4]{amsart}
\usepackage{amsmath}
\usepackage{amsthm}
\usepackage{amssymb}
\usepackage{graphicx}
\usepackage{multirow}
\usepackage{lscape}
\usepackage{slashbox}
\usepackage[all]{xy}

\newtheorem{rem}{Bemerkung}

\newcommand{\Q}{{\mathbb Q}}
\newcommand{\G}{{\mathbb G}}
\newcommand{\Z}{{\mathbb Z}}
\newcommand{\F}{{\mathbb F}}
\newcommand{\N}{{\mathbb N}}
\newcommand{\C}{{\mathbb C}}
\newcommand{\R}{{\mathbb R}}
\newcommand{\BP}{{\mathbb P}}

\newcommand{\MF}{{\mathcal F}}
\newcommand{\MA}{{\mathcal A}}
\newcommand{\MB}{{\mathcal B}}

\newcommand{\MO}{{\mathcal O}}
\newcommand{\MT}{{\mathcal T}}

\newcommand{\MW}{{\mathcal W}}
\newcommand{\MS}{{\mathcal S}}

\newcommand{\La}{{\Lambda}}
\newcommand{\De}{{\Delta}}

\newcommand{\la}{{\lambda}}
\newcommand{\OM}{{\Omega}}
\newcommand{\Ga}{{\Gamma}}
\newcommand{\ga}{{\gamma}}
\newcommand{\bo}{{\boldsymbol{\omega}}}

\newcommand{\bnu}{{\boldsymbol{\nu}}}
\newcommand{\bC}{{\boldsymbol{C}}}
\newcommand{\bF}{{\boldsymbol{F}}}
\newcommand{\bX}{{\boldsymbol{X}}}
\newcommand{\om}{{\omega}}

\newcommand{\ba}{{\boldsymbol{a}}}
 \newcommand{\bx}{{\boldsymbol{x}}}
\newcommand{\bs}{{\boldsymbol{s}}}
\newcommand{\bk}{{\boldsymbol{k}}}
\newcommand{\bso}{{\boldsymbol{o}}}
\newcommand{\bsi}{{\boldsymbol{i}}}
\newcommand{\bsj}{{\boldsymbol{j}}}
\newcommand{\bal}{{\boldsymbol{\alpha}}}
\newcommand{\bla}{{\boldsymbol{\lambda}}}
\newcommand{\bA}{{\boldsymbol{A}}}
\newcommand{\bB}{{\boldsymbol{B}}}
\newcommand{\bD}{{\boldsymbol{D}}}

\newcommand{\lra}{{\longrightarrow}}

\makeindex
\numberwithin{equation}{section}

\begin{document}
\setcounter{section}{-1}
\parindent=0pt

\title[On Drinfeld modular forms of higher rank II]
{On Drinfeld modular forms of higher rank II} 
 \author{Ernst-Ulrich Gekeler}
 \maketitle
 \begin{abstract}
We show that the absolute value $|f|$ of an invertible holomorphic function $f$ on the Drinfeld
symmetric space $\OM^r$ $(r \geq 2)$ is constant on fibers of the building map to the
Bruhat-Tits building $\MB\MT$. Its logarithm $\log|f|$ is an affine map on the realization
of $\MB\MT$. These results are used to study the vanishing loci of modular forms 
(coefficient forms, Eisenstein series, para-Eisenstein series) and to determine their
images in $\MB\MT$.
 \vspace{0.3cm}\end{abstract}

 \section{Introduction}

The present paper continues the work reported in \cite{11}. There we started with the
investigation of growth/decay properties of (Drinfeld) modular forms for the group 
$\Ga = {\rm GL}(r,\F_q[T])$, where $\F_q$ is the field with $q$ elements and $r \geq 3$. 
(The case $r = 2$ is well understood since the 1990s, see \cite{6}, \cite{7}, \cite{9}).
Several results not essentially needed for the reasoning of \cite{11} were stated there
without proof, e.g.: If $f$ is an invertible holomorphic function on the Drinfeld symmetric space 
$\OM^r$ on which $\Ga$ acts, then 
 \begin{itemize}
 \item $\log_q|f(\bo)|$ is constant on fibers of the building map
 $$\la:\: \OM^r \lra \MB\MT$$
to the Bruhat-Tits building $\MB\MT$;
 \item regarded as a function on the set $\MB\MT(\Q)$ of points of $\MB\MT$ with rational 
barycentric coordinates, this function is affine (that is, interpolates linearly in simplices).
 \end{itemize}
As these fundamental facts turn out to be crucial for subsequent work, we give here complete proofs:
see Theorems 2.4 and 2.6. We use these and some results about functional determinants to study
the vanishing loci $V(f)$ of several classes of modular forms $f$ for $\Ga$. We are able to determine
the image $\la(V(f))$ under $\la$ for 
 \begin{itemize}
 \item $f$ one of the coefficient forms $g_1,\ldots,g_{r-1}$ (which together with the discriminant
$\De = g_r$ generate the algebra of modular forms of type 0: Theorem 4.2;
 \item $f$ an Eisenstein series $E_k$: Theorem 4.5;
 \item $f$ a para-Eisenstein series $\alpha_k$: Theorem 4.8.
 \end{itemize} 
We finally show that (after some normalization) the $k$-th para-Eisenstein series $\alpha_k$ is
a locally uniform limit of coefficient forms $_a\ell_k$, when the degree of $a \in \F_q[T]$ tends to
infinity: Theorem 4.8.
 \vspace{0.3cm}

Here we take the opportunity to point to recent articles of Basson \cite{2} and Basson-Breuer \cite{3}
about higher rank Drinfeld modular forms, whose results are complementary to the present work. 
 \vspace{0.3cm}

The plan of the paper is as follows.
 \vspace{0.3cm}

In Section 1, apart from recalling basic facts and definitions, we use successive minimum bases of
lattices in $C_{\infty}$ (the characteristic-$p$ analogue of the complex numbers) to describe the 
fundamental domain $\bF$ for $\Ga$ and to define the spectrum ${\rm spec}(\La)$, a fundamental
invariant of the $\F_q$-lattice $\La$.
 \vspace{0.3cm}

In Section 2 we give a geometric description of the fiber $\la^{-1}(\bx)$ of 
$\bx \in \MB\MT(\mathbb Q)$.
It turns out to be an affinoid with the absolute value property AVP: If $f$ is an invertible function on
$\la^{-1}(\bx)$, then the absolute value $|f(\bo)|$ is constant: Theorem 2.4. The fact that
$\log_q|f|$ is an affine function on $\MB\MT(\mathbb Q)$ (Theorem 2.6) is obtained by investigating
the restriction of $f$ to $\la^{-1}(\bs)$, where $\bs$ is a line segment inside a simplex $\sigma$ of
$\MB\MT$. In this situation we may reduce the assertion to a known fact of one-dimensional geometry due to
Motzkin \cite{14}, see also \cite{4} I.8.3.
 \vspace{0.3cm}

In Section 3 we study functional determinants
 $$\det_{1\leq i,j\leq r'} (\frac{\partial f_i}{\partial \om_j}),\: 1 \leq r' \leq r-1,\: \om_1,\ldots,\om_{r-1}
 \mbox{ the coordinates on } \OM^r,$$
where the $f_i$ are either
 \begin{itemize}
 \item the para-Eisenstein series $\alpha_i = \alpha_i(\La)$, or 
 \item the Eisenstein series $E_{q^i-1}(\La)$, or 
 \item the coefficient forms $_a\ell_i$ for some fixed $a \in \F_q[T]$ (in particular, $_T\ell_i = g_i$, the
basic coefficient forms which describe the generic Drinfeld module $\phi^{\bo}$ of rank $r$).
 \end{itemize}
From considering the $T$-torsion $_T\phi^{\bo} \cong \F_q^r$ of $\phi^{\bo}$, we reduce the
calculation of the functional determinant to evaluating certain Moore determinants. This notably
shows the non-vanishing (Proposition 3.15).
 \vspace{0.3cm}

In the fourth section we apply the preceding to locate the vanishing sets $V(f)$ inside $\bF$ and their
images $\la(V(f))$.
 \vspace{0.3cm}

For each function $f$ of type $g_i$, $E_k$, $\alpha_i$, there is a natural and explicitly computable range
$R(f)$ such that $\la(V(f))$ is contained in $R(f)$. It is fairly easy to verify that in fact 
 \begin{itemize}
 \item $\la(V(g_i))=R(g_i) := \MW_{r-i}$ (Theorem 4.2), where $\MW_{r-1}$ is the $(r-1)$-th wall of
$\MW = \la(\bF)$ 
 \end{itemize}
and
 \begin{itemize}
 \item $\la(V(E_k)) = R(E_k) := \MW_{r-1}$ (Theorem 4.5; $E_k$ is the Eisenstein series of weight $k$ with
$0 < k \equiv 0 \,(\bmod q-1)$,
 \end{itemize}
but difficult to show that
 \begin{itemize}
 \item $\la(V(\alpha_i)) = R(\alpha_i) := \MW(i)$ (Theorem 4.8).
 \end{itemize}
Here $\MW(i) = \la(\MF(i))$ is defined through spectral properties of the corresponding lattices.
 \vspace{0.3cm}

These theorems also include smoothness and intersection properties of the $V(f)$. Together with
Theorems 2.4 and 2.6 and Remarks 2.9 and 4.15, they allow a precise description of $|f(\bo)|$ or
$\|f\|_{\bx}$ for any of these modular forms $f$, $\bo \in \bF$ and $\bx \in \MW(\mathbb Q)$, where
$\|f\|_{\bx}$ is the spectral norm of $f$ on $\la^{-1}(\bx)$.
 \vspace{0.5cm}\\
{\bf Notation.}
 \vspace{0.3cm}

We use essentially the same notation as in \cite{11}, that is:
 \vspace{0.3cm}

$\F=\F_q =$ finite field with $q$ elements, of characteristic $p$\\
$\overline{\F} =$ algebraic closure of $\F$, $\F^{(n)} = \{x \in \overline{\F}~|~x^{q^n}=x\}$\\
$A = \F[T]$ the polynomial ring over $\F$, $K= {\rm Quot}(A) = \F(T)$\\
$K_{\infty} = \F((T^{-1}))$ the completion of $K$ w.r.t. the absolute value $|~|$ at infinity, normalized
by $|T| = q$\\
$C_{\infty} =$ completed algebraic closure of $K_{\infty}$, $O_{\infty} \subset K_{\infty}$ and
$O_{C_{\infty}} \subset C_{\infty}$ the rings of integers\\
$\log = -v_{\infty}:\: C_{\infty}^{\ast} \lra \mathbb Q$ the map $z \longmapsto \log_q|z|$\\
$\OM^r = \{\bo = (\om_1:\cdots:\om_r) \in \BP^{r-1}(C_{\infty})~|~ \mbox{ the $\om_i$ or $K_{\infty}$-{\rm l.i.}}\}$,
where $r \geq 2$ and l.i. is short for linearly independent\\
$\OM^r(R) = \{\bo \in \OM^r~|~\mbox{ the $\om_i$ lie in the subring $R$ of $C_{\infty}$}\}$\\
$\Ga = {\rm GL}(r,A)$ the modular group with center $Z \cong \F^{\ast}$\\
$\Ga(T) = \{\gamma \in \Ga~|~ \gamma \equiv 1\,(\bmod\,T)\}$ the $T$-th congruence subgroup \\
$\tau$ is a non-commutative variable subject to $\tau c = c^q\tau$ for $c \in C_{\infty}$. 
We identify the $\F$-algebra $C_{\infty}\{\tau\}$ of ``polynomials'' in $\tau$ with the $\F$-algebra 
${\rm End}_{\F}(\G_a/C_{\infty}) = \{\underset{\rm finite}{\sum}a_iX^{q^i}~|~a_i \in C_{\infty}\}$ of $q$-additive 
polynomials via $\tau^i \leftrightarrow X^{q^i}$; ditto for ``power series'' $C_{\infty}\{\{\tau\}\} 
\stackrel{!}{=} \{\underset{i \geq 0}{\sum} a_iX^{q^i}\}$.
 \vspace{0.3cm}

Given an $\F$-lattice $\La$ in $C_{\infty}$,
 $$\begin{array}{lll}
{\displaystyle e_{\La}(X)}& = & {\displaystyle X \prod_{\la \in \La}{}^{'}(1-X/\la)= 
 \sum_{0\leq i \leq \dim_{\F}\La} \alpha_i(\La)X^{q^i} = \sum \alpha_i(\La) \tau^i} 
  \vspace{0.2cm}\\
\log_{\La}(X) & = & {\displaystyle \sum \beta_i(\La) \tau^i = \mbox{inverse of $e_{\La}(X)$ in $C_{\infty}\{\{\tau\}\}$}}
 \vspace{0.2cm}\\
E_k(\La) &=& {\displaystyle \sum_{\la\in \La}{}^{'}\la^{-k}}
 \end{array}$$
are the exponential function, the logarithm functions, the $k$-th Eisenstein series of $\La$, respectively.
Here, as usual, a primed product $\prod'$ or sum $\sum'$ means the product or sum over the non-vanishing
elements of the index set.
 \vspace{0.3cm}

If $\om_1,\ldots,\om_r$ are $K_{\infty}$-l.i. with lattice $\La = \La_{\bo} = \sum A\om = \sum A\om_i$ then 
$e_{\bo}:= e_{\La}$, and $\phi^{\bo} := \phi^{\La}$ denotes the attached Drinfeld module with operator polynomial
 $$\phi_a^{\bo}(X) = \sum_{0\leq k \leq r \cdot \deg\,a} {_a\ell_k}(\bo) X^{q^k} = \sum {_a\ell_k} (\bo) \tau^k$$
and $a$-torsion submodule $_a\phi^{\bo}$ ($a \in A$).
 \vspace{0.3cm}

Reduced analytic subspaces $\bX$ of $\BP^n(C_{\infty})$ or $\mathbb A^n(C_{\infty})$ are usually described through
their sets $\bX(C_{\infty})$ of $C_{\infty}$-valued points. Similarly, we often don't distinguish in notation
between a simplicial complex $\MS$ and its realization $\MS(\R)$.
 \vspace{0.3cm}

Finally, the cardinality of the finite set $X$ is denoted by $\#(X)$, the multiplicative group of the ring
$R$ by $R^{\ast}$.
 
 \section{Modular Forms}

(1.1) Recall that the Drinfeld symmetric space for $r \geq 2$ is 
 $$\begin{array}{lll} \OM^r &:=& \{\bo = (\om_1:\cdots: \om_r) \in \BP^{r-1}(C_{\infty})~|~ \om_1,\ldots,\om_r 
   \mbox{ l.i. over } K_{\infty}\}\\
 & & \mbox{(l.i. = linearly independent)} \\
 & = & \BP^{r-1}(C_{\infty}) \setminus \underset{H~ {\rm hyperplane~defined~over~ K_{\infty} }}
 {\bigcup\:H\,,}
 \end{array}$$
which carries a natural strucure as a rigid analytic space defined over $K_{\infty}$. The orbit
space $\Ga\setminus\OM^r$ is (the set of $C_{\infty}$-points of) the moduli space for Drinfeld
$A$-modules of rank $r$, as is explained below. Such a Drinfeld module $\phi$ is given through the 
operator polynomial
 $$\begin{array}{lll}
 \phi_T(X) &=& TX+g_1X^q+ \cdots + g_{r-1}X^{q^{r-1}}+\De X^{q^r}\\
 &=& T\tau^0+g_1\tau+ \cdots + g_{r-1}\tau^{r-1}+\De\tau^r,
 \end{array}\leqno{(1.1.1)}$$
where $\tau$ denotes the operator ($x \longmapsto x^q$) in ${\rm End}_{\F}(\G_a)$, with 
$g_i,\De \in C_{\infty}$ and $\De \not= 0$. We also put $g_0 := T$ and $g_r := \De$. Each such
$\phi = \phi^{\bo}$ comes from a uniquely determined $A$-lattice $\La = \La_{\bo} = A\om_1+ \cdots + A\om_r$
in $C_{\infty}$, where $\bo = (\om_1, \ldots,\om_r) \in C^r_{\infty}$ determines a point 
$(\om_1:\ldots:\om_r) \in \OM^r$. In this way, $A$-lattices of rank $r$ in $C_{\infty}$ (resp. homothety classes
of such lattices) correspond to rank-$r$ Drinfeld modules (resp. isomorphism classes of such modules). We
normalize projective coordinates on $\OM^r$ such that $\om_r = 1$. Then $\Ga$ acts on $\OM^r$ through
 $$\begin{array}{l}
 \gamma\bo = \bo',\: \om'_i = {\rm aut}(\ga,\bo)^{-1} {\displaystyle \sum_{1 \leq j \leq r}} \ga_{i,j}\om_j,
    \mbox{ where}\\
 {\rm aut}(\ga,\bo) = \ga_{r,1}\om_1+ \cdots + \ga_{r,r}\om_r \quad (\ga = (\ga_{i,j})),
 \end{array}\leqno{(1.1.2)}$$
and the $g_i = g_i(\bo)$ become functions on $\OM^r$ via
 $$\phi_T^{\bo} = \sum_{0 \leq i \leq r} g_i(\bo)\tau^i. \leqno{(1.1.3)}$$
In fact, $g_i$ is a modular form of weight $q^i-1$ and type 0 for the modular group $\Ga$.
 \vspace{0.3cm}

(1.2) A {\em modular form of weight} $k \in \N_0$ {\em and type} $m \in \Z/(q-1)$ for $\Ga$ is a 
function $f:\: \OM^r \lra C_{\infty}$
that
 \begin{itemize}
 \item[(i)] is holomorphic;
 \item[(ii)] satisfies
 $$f(\ga\bo) = \frac{{\rm aut}(\ga,\bo)^k}{({\rm det}\ga)^m} f(\bo) \quad \mbox{for } \ga \in \Ga$$
and
 \item[(iii)] satisfies a certain boundary condition (see (1.7), (1.8)).
 \end{itemize}
(Apart from the considerations of (1.7) and (1.8), modular forms of non-trivial types will play no role in this article.)
 \vspace{0.3cm}

(1.3) Let $\overline{M}^r$ be the weighted projective space ${\rm Proj}\,C_{\infty}[X_1,\ldots,X_r]$, where
the weight of $X_i$ is defined as $wt(X_i) := w_i := q^i-1$. Then the moduli scheme for rank-$r$ Drinfeld
modules is the open subscheme $M^r$ of $\overline{M}^r$ defined by the non-vanishing of $X_r$:
 $$M^r := ({\rm Proj}\,C_{\infty}[X_1,\ldots,X_r])_{X_r \not= 0}, $$
the map
 $$\begin{array}{lll}
 \OM^r & \lra & M^r(C_{\infty})\\
 \bo & \longmapsto & (g_1(\bo): \cdots : g_r(\bo))
 \end{array}$$
is analytic and $\Ga$-invariant, and defines an isomorphism of analytic spaces
 $$\Ga\setminus \OM^r \stackrel{\cong}{\lra} M^r(C_{\infty}).$$

(1.4) For each $\F$-lattice (= discrete $\F$-subspace) $\La$ in $C_{\infty}$, we write its exponential function
$e_{\La}$ and its inverse $\log_{\La}$ in the non-commutative ring $C_{\infty}\{\{\tau\}\}$ as 
 $$\begin{array}{lll}
 e_{\La} &=& {\displaystyle \sum_{i\geq0} \alpha_i\tau^i = \sum\alpha_iX^{q^i} = X\prod_{\la\in \La}{}^{'}
   (1-X/\la)}\vspace{0.2cm}\\
 \log_{\La} &=& \underset{i \geq 0}{\sum} \beta_i \tau^i.
 \end{array}$$
Then
 $$\alpha_0 = \beta_0 = 1,\:\sum_{i+j=k} \alpha_i\beta_j^{q^i} = \sum_{i+j=k} \alpha_i^{q^j}\beta_j=0 \mbox{ for }
  k>0,\leqno{(1.4.1)}$$
and $-\beta_j$ agrees with the Eisenstein series $E_{q^j-1}$, where 
$$E_k = \sum_{\la \in \La}{}^{'} \la^{-k} \: (k > 0) \mbox{ and } E_0 = -1.\leqno{(1.4.2)}$$
The quantities $\alpha_i,\beta_j,E_k$ depend on $\La$ and are written as $\alpha_i(\La),\ldots$ 
(or $\alpha_i(\bo),\ldots$ if $\La$ happens to be an $A$-lattice with $A$-basis $\{\om_1,\ldots,\om_r\}$).
Given such an $A$-lattice $\La$ and $a \in A$, we write the $a$-th operator polynomial of the associated
Drinfeld module $\phi=\phi^{\bo}$ as 
 $$\phi_a^{\bo} = \sum_{0 \leq i \leq r \cdot \deg\,a}{_a\ell_i}\tau^i,\leqno{(1.4.3)}$$
with ${_a\ell_i} = {_a\ell_i}(\La) = {_a\ell_i}(\bo),\: \bo = (\om_1,\ldots,\om_r)$. In particular,
${_T\ell_i} = g_i$. All the functions $\alpha_k$, $\beta_k$, $E_k$ are modular forms of type 0 for $\Ga$,
with weight $q^k-1$ for $\alpha_k$, $\beta_k$ and $k$ for $E_k$.
 \vspace{0.3cm}

(1.5) A {\em successive minimum basis} (SMB) of the $A$-lattice $\La$ of rank $r$ is an ordered 
(we (ab)use the curly brackets notation of non-ordered sets) $A$-basis $\{\om_1,\ldots\om_r\}$ which 
satisfies for $1 \leq i \leq r$:
 \vspace{0.3cm}

$|\om_i|$ is minimal among $\{|\la|~|~ \la \in \La \setminus A-\mbox{span of } \{\om_1,\ldots,\om_{i-1}\}\}$.
Such an SMB exists for each $A$-lattice $\La$, and it has the properties (\cite{10}, Proposition 3.1):
 \begin{itemize}
 \item[(i)] for $a_1,\ldots,a_r \in K_{\infty}$, $|\underset{1 \leq i \leq r}{\sum} a_i\om_i| = 
\underset{i}\max|a_i\om_i|$;
 \item[(ii)] the series $|\om_1|,\ldots,|\om_r|$ is an invariant of $\La$ and doesn't depend on the choice of
the SMB $\{\om_1,\ldots,\om_r\}$.
 \end{itemize}

(1.6) Let $\bF$ be the set
 $$\bF = \{\bo \in \OM^r~|~ \{\om_r,\ldots,\om_1\} \mbox{ is an SMB of } \La_{\bo} = \sum_{1 \leq i \leq r} A\om_i\}$$
(note the reverse order!). It is an admissible open subspace of the analytic space $\OM^r$, and each $\bo \in \OM^r$ 
is $\Ga$-equivalent with at least one and at most finitely many $\bo' \in \bF$. We call $\bF$ the {\em fundamental
domain} for $\Ga$ on $\OM^r$. As modular forms $f$ are uniquely determined on $\F$, we will focus our study to
the restriction of $f$ to $\bF$.
 \vspace{0.3cm}

{\bf Remark.} Of course, the defining condition for $\bo \in \bF$ doesn't depend on the choice of projective
coordinates for $\bo$. This notion of fundamental domain is weaker than the requirements on classical
fundamental domains, as the almost uniqueness of representatives $\bo' \in \bF$ cannot be achieved. Let
e.g. $\ga \in {\rm GL}(r,\bF) \hookrightarrow \Ga$ be an upper triangular matrix. Then with $\bo$ also
$\ga \bo$ belongs to $\bF$.
 \vspace{0.3cm}

(1.7) We may now specify the missing boundary condition (iii) in (1.2). Let $f:\: \OM^r \lra C_{\infty}$ be a function
that satisfies conditions (i) and (ii), with the type $m=0$ to fix ideas. We denote by $f^*$ the unique extension
of weight $k$ ($f^*(c\bo) = c^{-k}f^*(\bo)$ for $c \in C_{\infty}^*$) of $f$ to
 $$\OM^{r,*} := \{\bo = (\om_1,\ldots,\om_r) \in C_{\infty}^r~|~\mbox{the $\om_i$ are $K_{\infty}$-l.i.}\}.$$
Then $f^*$ is $\Ga$-invariant, and there exists a holomorphic function \linebreak
$F: \: \{(g_1,\ldots,g_r) \in C_{\infty}^r~|~g_r \not=0\} \lra C_{\infty}$ such that
 $$f^*(\bo) = F(g_1(\bo),\ldots,g_r(\bo)).\leqno{(1.7.1)}$$
The boundary condition (iii) now requires that $F$ admit a holomorphic extension to $C_{\infty}^r\setminus \{0\}$. 
Suppose this holds. Then by GAGA (or the non-archimedean Chow lemma, or by direct proof), $F$ is a polynomial
in the $g_i$, necessarily isobaric of weight $k$ (where $wt(g_i) = w_i = q^i-1$). Let such a polynomial $F$ be
given. Then the function
 $$\begin{array}{llll}
 f: & \OM^r & \lra & C_{\infty}\\
 & \bo & \longmapsto & F(g_1(\bo),\ldots,g_r(\bo)) \end{array}$$
is bounded on $\bF$, as the $g_i(\bo)$ are (\cite{11}, Corollary 4.16).
 \vspace{0.3cm}

Next, let $f$ with (i), (ii) be given and suppose it is bounded on $\bF$. Then $F:\: C_{\infty}^{r-1} \times
C_{\infty}^* \lra C_{\infty}$ defined by (1.7.1) extends to $C_{\infty}^r \setminus \{0\}$, due to the 
non-archimedean analogue of the Riemann removable singularities theorem \cite{1}.
 \vspace{0.3cm}

Now we allow the type $m$ to be non-trivial, $0 \le m < q-1$. Let $h$ be the function defined in
\cite{11}, Theorem 3.8 that satisfies
 $$(-1)^{r-1}h^{q-1} = T^{-1}g_r = T^{-1} \De.$$
It is a modular form of weight $w'_r := (q^r-1)/(q-1)$ and type 1. Some $f$ subject to (1.2)(i) and (ii) with $(k,m)$ 
arbitrary satisfies
 $$f^*(\bo) = F(g_1(\bo), \ldots,g_{r-1}(\bo), h(\bo))$$
with some holomorphic $F$, where $f^*(\ga \bo) = (\det\,\ga)^m f^*(\bo)$ for $\ga \in \Ga$. If $F$ extends
to $C_{\infty}^r \setminus \{0\}$, then $F$ is an isobaric polynomial of weight $k$, where $wt(h)=w'_r$.
The rest of the argument is as in the case where $m=0$. Therefore we have shown:
 \vspace{0.3cm}

{\bf 1.8 Proposition.} {\it Let $f:\: \OM^r \lra C_{\infty}$ be a function subject to conditions (i) and (ii) of
(1.2). The following are equivalent:
 \begin{itemize}
 \item[(a)] $f$ is regular along the divisor $X_r=0$ of $\overline{M}^r$ (that is, the associated $F$ extends);
 \item[(b)] $f$ is a polynomial $F$ in the forms $g_1,\ldots,g_{r-1},h$;
 \item[(c)] $f$ is bounded on the fundamental domain $\bF$. \hspace{2.8cm} $\Box$
 \end{itemize}
}

{\bf Remarks.} (i) The polynomial $F$ in (b) is necessarily isobaric of weight $k$ (weights $w_i$ for the $g_i$, weight
$w'_r$ for $h$), and of type $m$, that is $F(g_1,\ldots,g_{r-1},h) = h^mF'(g_1,\ldots,g_r)$, 
with some isobaric 
$F'$ of weight $k-m\cdot w'_r$.
 \medskip

(ii) While (a) and (b) are specific to the case considered (where the acting group is the full modular group
$\Ga = {\rm GL}(r,A)$ and the moduli scheme is easy to describe), condition (c) naturally generalizes.
If $\Ga'$ is some congruence subgroup of $\Ga$, modular forms $f$ for $\Ga'$ may be defined by the conditions (i),
(ii$\Ga'$), (iii$\Ga'$), where
 \begin{itemize}
 \item[(ii$\Ga'$)] the modular equation (ii) holds for $\ga \in \Ga'$;
 \item[(iii$\Ga'$)] $f$ is bounded on $\ga\bF$ for all $\ga$ in a system of representatives of $\Ga/\Ga'$
 \end{itemize}

This definition has the advantage that it doesn't require a precise description of the corresponding moduli
scheme.
 \vspace{0.3cm}

(1.9) Next, we consider  arbitrary $\F$-lattices $\La$ in $C_{\infty}$, that is, discrete 
(finite- or infinite-dimensional) $\F$-subspaces of $C_{\infty}$. A {\em successive minimum basis of}
$X$ {\em over} $\F$ (or $\F$-SMB for short) is an ordered $\F$-basis $\{\la_1,\la_2,\ldots\}$ with the
property analogous with (1.5): For each $i\in \N$ less than or equal to $\dim_{\F}(\La)$,
 $$|\la_i| \mbox{ is minimal among } \{|\la|~|~\la \in \La \setminus \F-\mbox{span of } \{\la_1,\ldots,\la_{i-1}\}\}.$$
As is easily seen, each $\La$ possesses an $\F$-SMB $\{\la_1,\la_2,\ldots\}$, and 
 \begin{itemize}
 \item[(i)] $|\sum a_i\la_i| = \underset{i}{\max}\{|\la_i|~|~ a_i \not= 0\}$ for $a_1,a_2,\ldots \in \F$, 
almost all vanishing;
 \item[(ii)] the series $|\la_1|, |\la_2|, \ldots$ is an invariant of $\La$ and independent of the choice
of the $\F$-SMB $\{\la_1,\la_2,\ldots\}$
 \end{itemize}
We call that series the {\em spectrum} ${\rm spec}(\La)$ of $\La$. Further, $\La$ is {\em separable} if
${\rm spec}(\La)$ is multiplicity free, i.e., $|\la_1| < |\la_2| < \ldots$, and {\em inseparable} otherwise,
$k$-inseparable if $|\la_k| = |\la_{k+1}|$. We put on record the observation:
 \vspace{0.3cm}

(1.10) Knowing ${\rm spec}(\La)$ is the same as knowing the Newton polygon $NP(e_{\La})$ of $e_{\La}$ 
(as defined in \cite{15} II Sect. 6), as by (i) both are equivalent to knowing the numbers of elements
of $\La$ of given sizes. In particular, $\La$ is separable if and only if the segments of $NP(e_{\La})$
have lengths $(q-1)$, $(q-1)q$, $(q-1)q^2$, ..., as in this case there are precisely $(q-1)q^{i-1}$
elements $\la \in \La$ with $|\la| = |\la_i|$.
 \vspace{0.3cm}

{\bf 1.11 Proposition.} {\it Let $\La$ be an $\F$-lattice with $\alpha_k(\La) = 0$ for some
$k < \dim_{\F}(\La)$. Then $\La$ is $k$-inseparable. Conversely, if $\La$ is $k$-inseparable, there
exists an isospectral $\F$-lattice $\La'$ (i.e., ${\rm spec}(\La) = {\rm spec}(\La'))$ such that 
$\alpha_k(\La') =0$.}
 \vspace{0.1cm}

 \begin{proof}
If $\alpha_k(\La) = 0$ then the abcissa $q^k$ cannot be a break point of $NP(e_{\La})$. Therefore,
$\#\{\la \in \La~|~|\la| = |\la_k|\} = \#\{\la \in \La~|~ |\la| = |\la_{k+1}|\}$, that is,
$|\la_k| = |\la_{k+1}|$. Let now $\La$ be $k$-inseparable, $e_{\La} = \underset{i \geq 0}{\sum}\alpha_i\tau^i$. 
Its Newton polygon doesn't change if we replace $\alpha_k$ with $\alpha'_k = 0$. The
lattice $\La' = {\rm ker}(e_{\La'})$ with $e_{\La'} = \underset{i\geq 0 \atop i\not= k}{\sum} 
\alpha_i\tau^i$ is as wanted.
 \end{proof}

{\bf Remark.} The same argument shows: If $\La$ is $k$-inseparable for all $k \in S$, where $S$ is a 
possibly infinite subset of $\N$, there exists an isospectral lattice $\La'$ with $\alpha_k(X') = 0$
for all $k \in S$.
 \vspace{0.3cm}

(1.12) Suppose we are given an $A$-lattice $\La$ with SMB $\{\om_r,\om_{r-1},\ldots,\om_1\}$. Then
$\{T^j\om_i~|~1 \leq i \leq r,\, j \in \N_0\}$ is an $\F$-basis of $\La$, from which we may construct
an $\F$-SMB by conveniently ordering the indices $(j,i)$:
 $$(j,i)\prec (j',i') \mbox{ if } |T^j\om_i| < |T^{j'}\om_{i'}| \mbox{ or } (|T^j\om_i| = T^{j'}\om_{i'}
  \mbox{ and } i > i').$$
Hence the $\F$-SMB starts
 $$\la_1 = \om_r,\: \la_2 = T\om_r, \ldots, \la_j = T^{j-1}\om_r,\, \la_{j+1} = \om_{r-1},$$
where $j$ is maximal such that $|T^{j-1}\om_r| \leq |\om_{r-1}|$.
 \vspace{0.3cm}

Recall that $\log\,z = \log_q|z|$ for $z \in C_{\infty|}^*$, so $\log\,T = 1$. In particular, $\La$ is separable if and
only if the $\log\,\om_i \in \Q$ are all incongruent modulo $\Z$ and $\La$ is 1-inseparable if and only if $|\om_{r-1}|
=|\om_r|$. For $r \geq 3$, $\La$ is 2-inseparable if 
 $$\begin{array}{cll}
 \mbox{$either$ } &  |\om_{r-2}| = |\om_{r-1}| < q|\om_r| & (\la_1 = \om_r,\, \la_2 = \om_{r-1},\,\la_3 = \om_{r-2})\\
 \mbox{$or$ } & |\om_{r-1}| = q|\om_r| & (\la_1 = \om_r,\,\la_2 = T\om_r,\, \la_3 = \om_{r-1}).
 \end{array}$$

 \section{A closer look to the building map}

(2.1) Let $\la:\: \OM^r \lra \MB\MT(\Q)$ be the building map as described in \cite{11} (2.3) onto the points
with rational barycentric coordinates of the Bruhat-Tits building $\MB\MT$ of ${\rm PGL}(r,K_{\infty})$. 
The apartment $\MA$ of $\MB\MT$ is the full subcomplex defined by the standard torus $T$ of diagonal matrices of 
${\rm GL}(r,K_{\infty})$, with set of vertices 
 $$\MA(\Z) = \{[L_{\bk}]~|~ \bk = (k_1,\ldots,k_r) \in \Z^r\},$$
where $[L_{\bk}]$ is the homothety class of the $O_{\infty}$-lattice 
$L_{\bk} = (T^{k_1}O_{\infty},\ldots,T^{k_r}O_{\infty})$ in $K_{\infty}^r$. We have $[L_{\bk}] = [L_{\bk'}] \Leftrightarrow
\bk'-\bk = (k,k,\ldots,k)$ for some $k \in \Z$. The realization $\MA(\R)$ (for which we henceforth briefly
write $\MA$) is an euclidean affine space with translation group
 $$(T(K_{\infty})/K_{\infty}^* T(O_{\infty})) \otimes \R \stackrel{\cong}{\lra} \R^r/\R(1,1,\ldots,1)
  \stackrel{\cong}{\lra} \{\bx \in \R^r~|~ x_r = 0\}\leqno{(2.1.1)}$$
and with the natural choice of origin $\bso= [L_{\boldsymbol 0}]$. That is, we use that isomorphism as a description of 
$\MA = \MA(\R)$. The choice of the Borel subgroup of upper triangular matrices in ${\rm GL}(r,K_{\infty})$ determines
the Weyl chamber
 $$\MW = \{\bx \in \MA~|~ x_i \geq x_{i+1} \mbox{ for } 1 \leq i < r\}\leqno{(2.1.2)}$$
with walls
 $$\MW_i = \{\bx \in \MW~|~x_i = x_{i+1}\} \quad (1 \leq i < r).$$
Then $\MW$ is a fundamental domain (in the classical sense) for $\Ga$ on $\MB\MT$, that is, each 
$\bx \in \MB\MT(\R)$
is $\Ga$-equivalent with a unique $\bx \in \MW$. The relationship between the fundamental domains $\bF$ in
$\OM^r$ and $\MW$ in $\MB\MT$ is 
 $$\la(\bF) = \MW(\Q),\: \la^{-1}(\MW) = \MF.\leqno{(2.1.3)}$$
We define
 $$\bF_i := \la^{-1}(\MW_i) = \{\bo \in \bF~|~|\om_i| = |\om_{i+1}|\}$$
and for $\bx \in \MW(\Q)$
 $$\bF_{\bx} := \la^{-1}(\bx) = \{\bo \in \bF~|~\log\,\om_i = x_i,\: 1 \leq i \leq r\}.$$
These are admissible open subspaces of $\bF$, and $\bF_{\bx}$ even affinoid (see (2.4)).
 \vspace{0.3cm}

{\bf 2.2 Definition.} {\it Let $\bX = {\rm Sp}(B)$ be an open affinoid subspace of some affine or projective
space over $C_{\infty}$. The {\em spectral norm} of $f \in B$ is $\|f\|_{\bX} := 
\sup\{|f(x)|~|~x \in \bX\}$. The space $\bX$ satisfies the {\em absolute value property} AVP if 
and only if each unit $f \in B^*$ has constant absolute value on $\bX$.}
 \vspace{0.3cm}

{\bf 2.3 Examples.} 
 \begin{itemize}
 \item[(0)] It is well-known that the $r$-dimensional unit ball 
 $$\{(\om_1,\ldots,\om_r) \in C_{\infty}^r~|~|\om_i| \leq 1 \mbox{ for all } i\}$$
satisfies AVP.
 \item[(i)] Let $\bC_a := \{\om \in C_{\infty}~|~ |\om| = q^a\}$ be the circumference with radius $q^a$, 
$a \in \Q$. Then $\bC_a$ satisfies AVP (\cite{12} p. 93).
 \item[(ii)] For $r \geq 2$ we put
 $$\hspace*{1.3cm}\bX^r := \{\bo = (\om_1: \cdots : \om_r) \in \BP^{r-1}(O_{C_{\infty}})~|~|\ell_H(\om_1,\ldots,\om_r)| < 1\},$$
where $H$ runs through the finite set of hyperplanes of $\BP^{r-1}(\F)$ and $\ell_H:\: \F^r \lra \F$ is a linear
form with kernel $H$, uniquely extended to an $O_{C_{\infty}}$-linear form $O_{C_{\infty}}^r \lra O_{C_{\infty}}$.
Then $\bX^r = \la^{-1}(\boldsymbol 0)$, where $\boldsymbol 0 = (0,\ldots,0)$ is the origin of 
$\MW(\Q) \subset \MA(\Q)$, and $\bX^r$ 
satisfies AVP (\cite{11} (2.5), (2.7)).
 \item [(iii)] If $\bX_1,\ldots,\bX_s$ satisfy AVP, so does $\bX := \bX_1 \times \cdots \times \bX_s$.
This is obvious, since a unit of $\bX$ gives rise to units on the $\bX_i$ by fixing the $j$-coordinates,
$j \not=i$.
 \item[(iv)] AVP is far from being satisfied for general affinoids. Let for example $\bX$ be the annulus
 $$\bA_{a,b} = \{\om \in C_{\infty}~|~a \leq \log\,\om \leq b\}$$
with rational numbers $a < b$. The coordinate function ``$\om$'' is a unit with non-constant absolute
value.
 \vspace{0.3cm} \end{itemize}

{\bf 2.4 Theorem.} {\it For each $\bx \in \MB\MT(\Q)$, the inverse image $\la^{-1}(\bx)$ is an open affinoid
subspace of $\OM^r$ that satisfies AVP.}
 
 \begin{proof}
(i) We may assume $\bx \in \MW(\Q)$ and thus $\la^{-1}(\bx) = \bF_{\bx} = \{\bo \in \bF~|~\log\,\om_i = x_i,\mbox{ all } i\}$ 
with $\bx = (x_1,x_2,\ldots,x_r)$, $x_1 \geq \cdots \geq x_r = 0$. Write $\{1,2,\ldots,r\}$ as the disjoint
union of classes $S$, where $i,j$ are in the same class if and only if $x_i \equiv x_j\,(\bmod\,\Z)$. For
each class $S$ let
 $$\bC_{\bx,S} := \{\bo_S = (\ldots,\om_i,\ldots)_{i\in S}~| \om_i \in C_{\infty},\, \log\,\om_i = x_i,\, \om_i=1
 \mbox{ if } i = r\}.$$
Then 
 $$\bF_{\bx} = \{\bo \in \prod_{S} \bC_{\bx,S}~|~ \bo \in \bF\},$$
where the condition ``$\bo \in \bF$'' is equivalent with 
 \begin{itemize}
 \item[(a)] $\bo \in \OM^r$ (i.e., $\om_1,\ldots,\om_r$ are $K_{\infty}$-linearly independent) and
 \item[(b)] $\{\om_r = 1,\, \om_{r-1}, \ldots,\om_1\}$ forms an SMB of $\La_{\bo} = 
\underset{1 \leq i \leq r}{\sum} A\om_i.$
 \end{itemize}
(ii) For each class $S$, consider the conditions on $\bo$:
 \begin{itemize}
 \item[(a$_S$)] the $\om_i$ ($i \in S$) are $K_{\infty}$-l.i. and
\item[(b$_S$)] the $\om_i$ ($i\in S$) in reverse order form an SMB of $\La_{\bo,S} := \underset{i\in S}
{\sum}  A\om_i$.
 \end{itemize}
It follows from (1.5)(i) and the construction of $\bC_{\bx,S}$ that (a) holds for $\bo$ whenever (a$_S$) holds
for all $S$.
 \medskip\\
Suppose that (a$_S$) and (b$_S$) hold for all $S$. By the above, $\om_1,\ldots,\om_r$ are $K_{\infty}$-l.i. and 
$|\om_r| \leq |\om_{r-1}| \leq \ldots |\om_1|$. We claim that, moreover, they form an SMB of $\La_{\bo}$.
 \medskip\\
For, suppose that there exists some $i < r$ and $\om \in \La_{\bo} \setminus 
\underset{i<j\leq r}{\sum} A\om_j$ with
$|\om_i| > |\om|$. Write
 $$ \om = \sum_{1 \leq j \leq r} a_j\om_j = \sum_{S} \om_S \mbox{ with } \om_S = \sum_{j\in S} a_j\om_j, \: a_j\in A.$$
Then $|\om_i| > |\om| = \underset{s}{\max} |\om_S|$.
 \medskip

If $i\in S$, then $\om_S \not=0$ contradicts (b$_S$). If $i \not\in S$ then $\om_S \in 
\underset{j\in S\atop 
|\om_j|<|\om_i|}{\sum}A\om_j \subset \underset{j\in S\atop j>i}{\sum} A\om_j$. Hence $\om = \sum \om_S \in \sum_{j>i} A\om_j$, which conflicts
with the assumption on $\om$.
 \medskip

(iii) By (i) and (ii) we have 
 $$\begin{array}{lll}
 \bF_{\bx} &=& \{\bo = (\ldots,\bo_S,\ldots) \in {\displaystyle \prod_{S}} \bC_{\bx,S}~|~\mbox{for each 
$S,\bo_S$ satisfies}\\
 & & \hspace*{5.6cm}\mbox{(a$_S$) and (b$_S$)}\}\\
 &=& {\displaystyle \prod_{S}} \bB_{\bx,S}, \end{array}$$
where
 $$\bB_{\bx,S} = \{\bo_S = (\om_i)_{i\in S} \in \bC_{\bx,S}~|~\mbox{(a$_S$) and (b$_S$) hold}\}.$$
We will see that all the $\bB_{\bx,S}$ are isomorphic with spaces treated in Example 2.3.
 \medskip

(iv) If \fbox{$S = \{r\}$} then $\bB_{\bx,S} = \{1\}$.
 \medskip\\
If \fbox{$r \in S \not=\{r\}$} then $\bB_{\bx,S} = \{\bo_S = (\om_i)_{i\in S}~|~\log\,\om_i = x_i,\mbox{ (a$_S$) and 
(b$_S$) hold}\}.$
 \medskip\\
As $x_i = \log\, \om_i \in \N_0$, we may scale $\om'_i := T^{-x_i} \om_i$ and find 
$$\bB_{\bx,S} \cong \{\bo_S = (\om_i)_{i\in S}~|~|\om_i|=1,\, \om_r = 1, 
 \mbox{ (a$_S$) and (b$_S$) hold}\},$$
which equals the space $\bX^s$ with $s:= \#(S)$ of Example 2.3(ii), for which AVP is satisfied.
 \medskip

If \fbox{$S = \{i\}$ with $i \not= r$} then (a$_S$) and (b$_S$) are trivially fulfilled, so
$\bB_{\bx,S} = \bC_{\bx,S} = \{\bo \in C_{\infty}~|~ \log\,\om_i = x_i\} = \bC_{x_i}$, 
a circumference, for which AVP holds by Example 2.3(i).
 \medskip

Let now $S$ be such that \fbox{$r \not\in S$ and $s:= \#(S)>1$}. As the $x_i=\log\,\om_i$ for
$i \in S$ are all congruent modulo $\Z$, we may scale the $\om_i$ by integral powers of $T$ such
that they have the same absolute value. Again multiplying with a fractional power of $T$, we may
achieve $|\om_i|=1$ for $i\in S$. Hence 
 $$\bB_{\bx,S}\cong \{\bo_S = (\om_i)_{i\in S}~|~|\om_i|=1, \mbox{ (a$_S$) and (b$_S$) hold}\},$$
which, by projecting to the last coordinate, becomes isomorphic with $\bX^s \times \boldsymbol C_0$
(see again Examples 2.3).
 \medskip

(v) Together, the analytic space $\bF_{\bx}$ is isomorphic with the product of affinoids each 
enjoying AVP, and thus enjoys the same properties.
 \end{proof}

{\bf Remark.} Let $\bx = (x_i)_{1 \leq i \le r}$ be an element of $\MW(\Q)$. Then $\bx$ lies in the
interior of a simplex $\sigma$ of maximal dimension $r-1$ if and only if all the $x_i$ are
incongruent modulo $\Z$, if and only if $\bF_{\bx}$ is a product of $(r-1)$ circumferences. 
Congruences $(\bmod\,\Z)$ of the $x_i$ imply that $\bx$ belongs to a simplex of smaller
dimension.
 \vspace{0.3cm}

{\bf 2.5 Corollary.} {\it If $f \in \MO(\OM^r)^*$ is a global unit (i.e., an invertible
holomorphic function on $\OM^r$) then $|f|$ factors through $\la$, and may thus be considered
as a function on $\MB\MT(\Q)$.} \hspace{5cm} $\Box$
 \vspace{0.3cm}

This holds e.g. for the discriminant function $\De$ and its root $h$. The values of $|\De|$ on $\MW(\Z)$ 
have been determined in \cite{11}. Since $\MW$ is a fundamental domain for $\Ga$ and $\De$ is modular,
this suffices to find $|\De(\bo)|$ for any $\bo$, in view of the next result.
 \medskip

For a unit $f$ we write $\log\,f(\bx)$ for the common value $\log_q|f(\bo)|$ of $\bo \in \bF_{\bx}$.
 \vspace{0.3cm}

{\bf 2.6 Theorem.} {\it Let $f$ be an invertible holomorphic function on $\OM^r$. Then
$\bx \longmapsto \log\,f(\bx)$ is an affine function on $\MB\MT(\Q)$. That is, given a
simplex $\sigma = \{\bk_0,\ldots,\bk_n\}$ of $\MB\MT$ and barycentric coordinates
$(t_0,\ldots,t_n) \in \Q_{\geq 0}^{n+1}$ with $\sum t_i = 1$ for
$\bx = \underset{0 \leq i \leq n}{\sum} t_i\bk_i$, then
 $$\log\,f(\bx) = \sum_{0 \leq i \leq n}t_i\log\,f(\bk_i).\leqno{(2.6.1)}$$
}

Before proving the theorem, we need some auxiliary results.
 \vspace{0.3cm}

{\bf 2.7 Proposition.} {\it Given rational numbers $a< b$, let $x_1,\ldots,x_s \in \bC_a$
and $y_1,\ldots,y_t \in \bC_b$ be finitely many points. For $1 \leq i \leq s$
(resp. $1 \leq j \leq t$) let
 $$\bB_i := \{z \in C_{\infty}~|~ |z-x_i| < q^{a_i}\}$$
be an open disc around $x_i$ with $a_i \leq a$ (resp. $\bB'_j$ an open disc of radius $q^{b_j}$
with $b_j \leq b$ around $y_j$). Suppose that all the $\bB_i,\bB'_j$ are disjoint. Let 
$f$ be an invertible holomorphic function on
 $$\bB := \bA_{a,b} \setminus(\bigcup_i \bB_i ~\cup ~\bigcup_j \bB'_j) = 
 \{\om \in C_{\infty}~|~q^a\leq |\om| \leq q^b,\,|\om-x_i| \geq a_i,\,|\om-y_j| \geq b_j\}.$$
Then there exists an invertible holomorphic function $g$ on $\bB$ with constant absolute
value and integers $m,m_1,\ldots,m_s,n_1,\ldots,n_t$ such that
 $$f(\om) = \om^m \prod_{i} (\om-x_i)^{m_i} \prod_{j} (\om-y_j)^{n_j} g(\bo).$$
}
 
\begin{proof}
This is Theorem I.8.5 of \cite{4}, adapted to our framework.
 \medskip\end{proof}

{\bf 2.8 Corollary.} {\it In the above situation, the absolute value $|f(\om)|$ depends only
on $|\om|$. Regarded as a function of $x:= \log\om$, $|\log\,f(\om)|$ interpolates linearly
between $a$ and $b$.}
 
 \begin{proof}
This results from the fact that $\om^m$, $(\om-x_i)^{m_i}$, $(\om-y_j)^{n_j}$ share these
properties.
 \medskip
\end{proof}

{\bf 2.9 Remark.} For further use, we note the following slight generalization. We suppress
the condition of $f$ to be invertible, but allow a finite number of zeroes in $\bC_a$ and
$\bC_b$. Replacing $|f(\om)|$ with the spectral norm $\|f\|_x$ on $\bC_x$, we still have 
 $$x \longmapsto \log_q\|f\|_x \mbox{ interpolates linearly between $a$ and $b$}.$$
This results from applying 2.8 to $\bB' = \bB$ minus the union of small disjoint open discs
around the zeroes of $f$. Note that $\log_q\|f\|_x = \log\,f(\om)$ for $a<x = \log\,\om < b$.
 
 \begin{proof} of Theorem 2.6. (i) We assume without restriction that $\sigma$ is the standard
simplex of maximal dimension $r-1$, with vertices $\boldsymbol 0 = (0,\ldots,0)$ and $\bk_i =
(1,1,\ldots 1,0,\ldots 0)$ with $i$ ones and $(r-i)$ zeroes ($1 \leq i < r$; we remind the reader
that all points $\bx = (x_1,\ldots,x_{r-1},0)$ have a redundant 0 as last entry $x_r$). Then 
 $$\la^{-1}(\sigma) = \{\bo \in \bF~|~ q \geq |\om_1| \geq |\om_2| \geq \cdots \geq |\om_{r-1}|
 \geq 1,\, \om_r = 1\}.$$ 
Recall that $\bo \in C_{\infty}^r$ gives rise to an element of $\bF$ if and only if (compare proof of
Theorem 2.4):
 \begin{itemize}
 \item[(a)] the entries $\om_i$ are $K_{\infty}$-l.i.;
 \item[(b)] $\{\om_r,\om_{r-1},\ldots,\om_1\}$ form an SMB of $\La_{\bo}$.
 \end{itemize}
(ii) Instead of $\bo = (\om_1,\ldots,\om_{r-1})$ we use $\bal =(\alpha_1,\ldots,\alpha_{r-1})$
as coordinates, where
 $$\alpha_i = \om_i/\om_{i+1} \mbox{ and so } \om_i = \prod_{i \leq j < r} \alpha_j.\leqno{(1)}$$
Then $|\alpha_i| \geq 1$ and $|\underset{1 \leq i < r}{\prod} \alpha_i| \leq q$. Fix some
$\bo^{(o)} = (\om_1^{(o)},\ldots,\om_{r-1}^{(o)}) \in \la^{-1}(\sigma)$ and for some $i$, $1 \leq i < r$, 
its $\alpha$-coordinates $(\alpha_1^{(o)},\ldots \alpha_{i-1}^{(o)}, \alpha_{i+1}^{(o)}, \ldots
\alpha_{r-1}^{(o)})$, where we suppose that 
 $$\tilde{q}_i := q/|\prod_{1\leq j < r\atop j\not= i} \alpha_j^{(o)}| > 1.$$
Consider the set of possible $\alpha_i$:
 $$\begin{array}{lll}
 \bA'_i &:=& \{\alpha \in \C_{\infty}~|~1 \leq |\alpha| \leq \tilde{q}_i\} \mbox{ and}\vspace{0.2cm}\\
 \bA_i &:= & \{\bo \in \bF~|~ \om_j/\om_{j+1} = \alpha_j^{(o)} \mbox{ for } j \not= i,\, \om_i/\om_{i+1} \in \bA'_i\},
 \end{array}$$
which, of course, depend on the choice of $\bo^{(o)}$.
 \medskip

(iii) The map
 $$\begin{array}{llll}
 u_i: & \bA_i & \lra & \bA'_i \\
 & \bo &\longmapsto & \alpha = \om_i/\om_{i+1} 
 \end{array}$$
is well-defined and injective. {\em Claim:} $u_i$ is an open embedding onto $\bA'_i$ minus a finite number of
disjoint open discs of shape
 $$\{\alpha \in C_{\infty}~|~|\alpha-z_j| < r_j\} \mbox{ with } |z_j| = 1 = r_j \mbox{ or } |z_j| = \tilde{q_i} = r_j.$$
To prove the claim, we consider for each $\alpha\in \bA'_i$ the point $\bo$ corresponding via (1) to
$(\alpha_1^{(o)}, \ldots, \alpha_{i-1}^{(o)}, \alpha, \alpha_{i+1}^{(o)}, \ldots, \alpha_{r-1}^{(o)})$. Then
$\alpha \in {\rm im}(u_i) \Leftrightarrow \bo \in \bF \Leftrightarrow$ conditions (a) and (b) hold for $\bo$.
 \medskip

First, suppose that \fbox{$1 < |\alpha| < \tilde{q}_i$}. Then $|\om_{i-1}^{(o)}| > |\om_i| > |\om_{i+1}^{(o)}|$, and
it is easily seen that (a) and (b) turn over from $\bo^{(o)}$ to $\bo$, hence $\alpha \in {\rm im}(u_i)$.
 \medskip

Next, assume \fbox{$|\alpha| = 1$}. Then $|\om_{i-1}^{(o)}| > |\om_i| = |\om_{i+1}^{(o)}|$. As 
 $$\begin{array}{lll}
\om_j &=& \om_j^{(o)} \mbox{ for } j > i\\
 \om_j &=& \om_j^{(o)} (\om_i/\om_i^{(o)}) \mbox{ for } j < i,
 \end{array}$$
in particular, $|\om_j| < q$ for all $j$, $\bo$ fulfills (a) and (b) if and only if $\om_i,\om_{i+1}^{(o)}, \ldots,
\om_k^{(o)}$ are $K_{\infty}$-linearly independent and form an SMB of their $A$-span, where $k$ with $i \leq k \leq r$ is maximal
with $|\om_i| = |\om_k^{(o)}|$. This is the case if and only if $|\om_i-\om_j^{(o)}| \geq |\om_i|$ 
for $j = i+1,\ldots,k$,
if and only if $|\alpha-\om_j^{(o)}/\om_{i+1}^{(o)}| \geq 1$ for $j = i+1,\ldots,k$.
 \medskip

If \fbox{$|\alpha| = \tilde{q}_i$}, then $|\om_{i-1}^{(o)}| = |\om_i| > |\om_{i+1}^{(o)}|$, and an analogous argument
shows that $\alpha \in {\rm im}(u_i)$ if and only if $|\alpha-\om_j^{(o)}/\om_{i+1}^{(o)}| \geq \tilde{q}_i$ for 
$j = k,\ldots,i-1$, where $k$ with $1 \leq k < i$ is minimal such that $|\om_k^{(o)}| = |\om_i|$.
 \medskip

Hence the claim is proved and $\bA_i$ is a one-dimensional analytic space isomorphic with some space 
that appears in
Proposition 2.7.
 \medskip

(iv) Consider the image $\la(\bA_i)$ in $\MB\MT(\Q)$. It is a maximal line segment in $\sigma$ parallel with
the vector $\bk_i = (1,1,\ldots,1,0\ldots0)$, and doesn't depend on $\bo^{(o)}$ itself, but only on the absolute 
values $|\om_i^{(o)}|$, or, what is the same, on $\la(\bo^{(o)})$.
 \medskip

Given our data $\sigma$ and $f$, an invertible function on $\la^{-1}(\sigma)$, a line segment $\bs$ in $\sigma$, 
(that is, in $\sigma(\Q)$, which we always implicitly assume) is called {\em well-behaved} if $\bx \longmapsto \log\,f(\bx)
:= \log_q\|f\|_{\bx}$ is a linear function along $\bs$. Applying Corollary 2.8 to $\bA_i$, we find that
$\la(\bA_i)$ is well-behaved. 
As the choices of $\bo^{(o)}$ and $i$ were arbitrary (the excluded case $\tilde{q}_i=1$
would lead to $\la(\bA_i) = \la(\bo^{(o)})$), this shows that each line $\bs$ in the simplex $\sigma$ 
which is parallel
with one of the base vectors $\bk_1,\ldots,\bk_{r-1}$ is well-behaved.
 \medskip

(v) With similar arguments one shows that each line segment $\bs$ in $\sigma$ parallel with a 1-simplex
 facing $\sigma$
and different from $({\boldsymbol 0},\bk_i)$ is well-behaved. The 1-simplices are of shape $\{\bk_i,\bk_j\}$ with 
$1 \leq i < j < r$, and one has to construct subspaces $\bA_{i,j}$ similar to the $\bA_i$ of (ii) such that
$\la(\bA_{i,j})$ equals a maximal line segment in $\sigma$ parallel with $\{\bk_i,\bk_j\}$. Alternatively,
one could use the fact that ${\rm PGL}(r,K_{\infty})$ acts transitively on the set of 1-simplices of $\MB\MT$. We omit
the details.
 \medskip

(vi) Now formula (2.6.1) holds along the 1-faces of $\sigma$, that is, along the boundary of each 2-face.
 \medskip

Let $\tau$ be an $n$-face of $\sigma$ ($n \geq 2$), and suppose that 
 \begin{itemize}
 \item[(a)] lines in $\tau(\Q)$ parallel with one of its 1-faces are well-behaved;
 \item[(b)] formula (2.6.1) holds on the boundary of $\tau$.
 \end{itemize}
Then an easy calculation shows that (2.6.1) holds for all points $\bx \in \tau (\Q)$. Therefore,
by induction on the dimension of faces of $\sigma$, we find that (2.6.1) holds for each $\bx \in
\sigma(\Q)$. This finishes the proof.
 \end{proof}

 \section{Functional determinants}

In this section, we deal with certain functional determinants related to the following families of
modular forms:
 \vspace{0.3cm}

(3.1) (i) the coefficient forms $g_k(\bo)$ ($1 \leq k \leq r$), see (1.1);
 \vspace{0.3cm}

(ii) the Eisenstein series
 $$\begin{array}{cll}
 E_k(\bo) &=& \underset{\la \in \La_{\bo}}{\sum{}^{'}}\la^{-k}, \mbox{ notably }\vspace{0.2cm}\\
 E_{q^k-1}(\bo) &=& -\beta_k(\bo), \mbox{ see (1.4)}; \end{array}$$

(iii) the para-Eisenstein series $\alpha_k(\bo)$, see (1.4);
 \vspace{0.3cm}

(iv) the coefficient forms $_a\ell_k(\bo)$ (see (1.4); in particular $_T\ell_k = g_k$);
\vspace{0.3cm}

(v) the forms $\mu_i(\bo) := e_{\bo} (\om_i/T)$ ($1\leq i \leq r$, see \cite{11} Sect. 3); they form an
$\F$-basis of the space $V:= {_T\phi^{\bo}}$ of $T$-division points of $\phi^{\bo}$; their reciprocals are
modular of weight 1 for the congruence subgroup $\Ga(T) = \{\ga\in \Ga~|~\ga \equiv 1\,(\bmod\,T)\}$.
In particular
 $$e_V(X) = T^{-1} \phi_T^{\bo}(X) = X+T^{-1} \sum_{1\leq j \leq r} g_j(\bo) X^{q^j},\leqno{(3.2)}$$
so the $g_j$ are rational functions of the $\mu_i$. In fact,
 \vspace{0.3cm}

(3.3) $T^{-1}g_j$ is the $(q^j-1)$-th elementary symmetric function $s_{q^j-1} \{\mu^{-1}\}$ of the set
$\{\mu^{-1}~|~0 \not= \mu \in V\}$, $V =$ set of $\F$-linear combinations of the $\mu_i$ ($1 \leq i \leq r$).
 \vspace{0.3cm}

These forms are subject to various relations; besides (3.3) we will use the following (see e.g. \cite{5},
Sect. 2; here and in the sequel, $\sum_{i+j=k} \ldots$ is short for $\underset{i,j\geq 0,i+j=k}{\sum} \ldots$):
 \vspace{0.3cm}

(3.4) (i)
 $$\begin{array}{rll}
 {\displaystyle \sum_{i+j=k} \alpha_i^{q^j} E_{q^j-1}}&=& 0 = {\displaystyle \sum_{i+j=k} \alpha_iE_{q^j-1}^{q^i},\:
   k>0,} \mbox{ so}\vspace{0.2cm}\\
 \alpha_k &=& {\displaystyle \sum_{0 \leq i<k} \alpha_i E_{q^{k-i}-1}^{q^i}};\end{array}
 $$
(ii)
 $$\begin{array}{ll}
 {\displaystyle \sum_{i+j=k} {_a\ell_i} \alpha_j^{q^i} = a^{q^k} \alpha_k,} \mbox{ that is}
  \vspace{0.2cm}\\
 {\displaystyle \sum_{0 \leq j < k} {_a\ell_{k-j}} \alpha_j^{q^{k-j}} = [a,k]\alpha_k,} \mbox{ where } 
 [a,k] := a^{q^k} -a. \end{array}$$

(3.5) We let $N^r$ be the ($C_{\infty}$-analytic space associated with the) moduli scheme for rank-$r$
Drinfeld modules $\phi$ with a structure of level $T$ (i.e., the choice of an ordered $\F$-basis
$\{\mu_1,\ldots,\mu_r\}$ of the space of $T$-division points $_T\phi$). As we can choose the $\mu_i$
arbitrary subject only to the condition of $\F$-linear independence,
 $$N^r = \BP^{r-1}(C_{\infty}) \setminus \bigcup\,H,$$
where $H$ runs through the finitely many hyperplanes defined over $\F$. Consider the commutative diagram
 $$\begin{array}{llcllll}
 \OM^r &\lra & \Ga(T)\setminus \OM^r & \stackrel{\cong}{\lra} & N^r &=& \{(\mu_1: \cdots : \mu_r) 
   \in \BP^{r-1} (C_{\infty}),\mu_1,\ldots,\mu_r\: \F-{\rm l.i.}\}  
   \vspace{0.3cm}\\
 & \searrow & \downarrow &&&& \vspace{0.3cm}\\
 & & \Ga\setminus \OM^r & \stackrel{\cong}{\lra} & M^r &=&({\rm Proj}(C_{\infty}[g_1,\ldots,g_r])_{g_r\not= 0},
 \end{array}\leqno{(3.6)}$$
where the maps are the natural quotient maps, in coordinates
 $$\begin{array}{lll}
 \bo & \longmapsto & (\mu_1(\bo): \cdots : \mu_r(\bo))
   \vspace{0.3cm}\\
 & \searrow & \downarrow 
  \vspace{0.3cm}\\
 & & (g_1(\bo): \cdots : g_r(\bo)).
 \end{array}$$
The right hand map is given by (3.3), and the group of the Galois covering on the right hand side is
 $$\Ga/\Ga(T) \cdot Z \stackrel{\cong}{\lra} {\rm GL}(r,\F)/Z = {\rm PGL}(r,\F),$$
which acts as matrix group on the column vector $(\mu_1,\ldots,\mu_r)^T$.

 \medskip

{\bf 3.7 Lemma.} {\it The quotient map $\OM^r \lra N^r$ is \'etale .}
 
 \begin{proof}
It suffices to show that $\Ga(T)$ acts without fixed points on $\OM^r$. So let $\ga \in \Ga(T)$ be such
that $\ga\bo = \bo$ for some $\bo \in \OM^r$. As $\ga$ fixes the image $\la(\bo) \in \MB\MT(\Q)$, $\ga$ has
finite order. Since $\ga \in \Ga(T)$, its eigenvalues all equal 1, so $\ga$ is unipotent. But a unipotent
$\ga \not= 1$ cannot have fixed points on $\OM^r$ (its eigenvectors on $C_{\infty}^r$ are $K$-rational), so $\ga =1$.
 \end{proof}

Note that any modular form $f$ for $\Ga$ is a rational function of $\mu_1,\ldots,\mu_r$, as follows from
Proposition 1.8 and (3.3).
 \vspace{0.3cm}

So $\frac{\partial f}{\partial \mu_j}$ is meaningful. In what follows, we regard ``forms'' as homogeneous functions
on the canonical cones
 $$\begin{array}{llll}
 \OM^{r,*} &:=& \{(\om_1,\ldots,\om_r) \in C_{\infty}^r~|~ (\om_1: \cdots : \om_r) \in \OM^r\} & \mbox{ over } \OM^r,
  \vspace{0.3cm}\\
 N^{r,*} &:=& \{(\mu_1,\ldots,\mu_r) \in C_{\infty}^r ~|~ (\mu_1:\cdots ,\mu_r) \in N^r\} & \mbox{ over } N^r,
  \end{array}$$
respectively. Then also $\frac{\partial f}{\partial \om_r}$ is meaningful.
 \vspace{0.3cm}

(3.8) Let now $D$ be one of the differential operators $\frac{\partial}{\partial \om_j}$ or 
$\frac{\partial}{\partial \mu_j}$ ($1 \leq j \leq r$). Then for $k > 0$
 $$D(\alpha_k)  = \sum_{0 \leq i < k} (E_{q^{k-i}-1})^{q^i} D(\alpha_i) + D(E_{q^k-1}),$$
as results from (3.4)(i). If $\bD(f)$ denotes either the vector 
$(\frac{\partial}{\partial \om_1}(f),\ldots, \frac{\partial}{\partial \om_r}(f))$ or the vector
$(\frac{\partial}{\partial \mu_1}(f), \ldots, \frac{\partial}{\partial \mu_r}(f))$ then
 \vspace{0.3cm}

(3.9) $\bD(E_{q^k-1}) = \bD(\alpha_k)+$ a linear combination of the $\bD(\alpha_i)$ with $1 \leq i < k$.
 \vspace{0.3cm}

In particular, we find
 \vspace{0.3cm}

{\bf 3.10 Proposition.} {\it For each natural number $r' \leq r$ and each $a \in A$ of degree $>0$, the three
functional determinants agree, where the indices $i,j$ always range between $1$ and $r'$:
 $$\begin{array}{lll}
 \det(\frac{\partial}{\partial \om_j}(E_{q^i-1})) &=& \det(\frac{\partial}{\partial \om_j}(\alpha_i))
  \vspace{0.2cm}\\
 &=& ([a,r'][a,r'-1] \cdots [a,1])^{-1} \det(\frac{\partial}{\partial\om_j}(_a\ell_i)).
 \end{array}$$
The statement remains true if all the $(\frac{\partial}{\partial \om_j})$ are replaced with 
$(\frac{\partial}{\partial \mu_j})$.}

 \begin{proof}
The first equality comes from (3.9). The second equality follows similarly, starting with (3.4)(ii), which leads
to 
 $$ \begin{array}{lll}
 [a,k]\bD(\alpha_k) &= &\bD(_a\ell_k)+ \mbox{ a linear combination of the } \bD(_a\ell_i) \\
 && \mbox{with }
 1 \leq i < k.\end{array}$$
The transfer - replacing the $\frac{\partial}{\partial \om_j}$ with the $\frac{\partial}{\partial \mu_j}$ - is
obvious.
 \end{proof}
 
Note that we can check the (non-) vanishing of these determinants with $a=T$, in which case $_a\ell_k = g_k$.
 \vspace{0.3cm}

(3.11) Next, we let $V$ be a finite $\F$-lattice in $C_{\infty}$ with an ordered $\F$-basis $\{\la_1,\ldots,\la_n\}$. All
the quantities $\alpha_k = \alpha_k(V) = \alpha_k(\la_1,\ldots,\la_n) = \alpha_k(\bla), \beta_k = \cdots = \beta_k(\bla) =
-E_{q^k-1} (\bla)$, $D_j := \frac{\partial}{\partial \la_j}$, $\bD(f) = (D_1(f),\ldots,D_n(f))$ refer to $V$.
From the identities in (1.4), we find for $n' \leq n$:
 $$\underset{1\leq i,j\leq n'}{\det} (D_j(\alpha_i)) = \underset{1 \leq i,j \leq n'}{\det} (D_j(E_{q^i-1})).
 \leqno{(3.12)}$$
We now recycle and generalize the argument of \cite{11} (6.3). 
 \vspace{0.3cm}

For an $\F$-linear map $\varphi:\: V \lra \F$ define 
 $$M(\varphi) := \sum_{\la \in V}{}^{'} \varphi(\la)/\la.$$
Then 
 $$D_j(E_{q^i-1}) = \sum_{\ba \in \F^n}{}^{'} \frac{a_j}{(a_1\la_1+ \cdots + a_n\la_n)} q^i = M(\varphi_j)^{q^i},$$
where
 $$\varphi_j: \, (a_1\la_1+ \cdots + a_n\la_n) \longmapsto a_j.$$
Hence the determinant in (3.12) is the Moore determinant (\cite{13} 1.9)
 $$\underset{1\leq i,j \leq n'}{\det} D_j(E_{q^i-1}(\bla)) = \underset{1\leq i,j \leq n'}{\det} 
 (M(\varphi_j)^{q^i}),$$
which doesn't vanish if and only if the $M(\varphi_j)$ ($1 \leq j \leq n'$) are $\F$-linearly independent. The latter holds
true for $1 \leq j \leq n$, as is shown by the argument (\cite{11}, Lemma 6.4). Therefore:
 \vspace{0.3cm}

{\bf 3.13 Proposition.} {\it The determinants appearing in $(3.12)$ never vanish.} \hspace{11.5cm}$\Box$
 \vspace{0.3cm}

{\bf Remark.} As the order of $\la_1,\ldots,\la_n$ is arbitrary, the non-vanishing of the functional determinants holds for
any given order.
 \vspace{0.3cm}

Let now $V$ be the $r$-dimensional $\F$-space $_T\phi^{\bo}$ with $\bo \in \OM^{r,*}$. Then by (3.2), $\alpha_i(V) = T^{-1} g_i(\bo)$, and
the $\mu_i(\bo)$ form an $\F$-basis of $V$. Combining (3.10) with (3.13) we find:
 \vspace{0.3cm}

{\bf 3.14 Proposition.} {\it Neither of the following determinants vanishes:
 $$\det(\frac{\partial}{\partial \mu_j}(\alpha_i)),\: \det(\frac{\partial}{\partial \mu_j}(E_{q^i-1})),\: 
\det(\frac{\partial}{\partial \mu_j} (_a\ell_i)).$$
Here ``$\det$'' $= \underset{1\leq i,j \leq r'}{\det}$, and $1 \leq r' \leq r$.}\hspace{5.4cm} $\Box$
 \vspace{0.3cm}

Now we return to our usual convention: modular forms are functions on $\OM^r$, whose coordinates are normalized by
$\om_r=1$. Then we are interested in functional determinants of rank $r' = r-1$. Fixing $i_0 \in \{1,2,\ldots,r\}$, the
$\mu_i$ ($1 \leq i \leq r$, $i \not= i_0$) are coordinates on $N^r$. From (3.7) we find that $\det(\frac{\partial\mu_i}
{\partial \om_j})$ never vanishes, where $1 \leq j \leq r-1$ and $1 \leq i \leq r$, $i \not= i_0$. Together with (3.10) and
(3.13) this yields:
 \vspace{0.3cm}

{\bf 3.15 Proposition.} {\it Let $a \in A$ be non-constant. Then the three functional determinants (where 
$1 \leq i,j \leq r-1$) never vanish:
 $$\det(\frac{\partial}{\partial \om_j} (E_{q^i-1})),\: \det(\frac{\partial}{\partial \om_j} (\alpha_i)),\: 
 \det(\frac{\partial}{\partial \om_j}(_a\ell_i)).$$}
\hspace*{12.3cm} $\Box$

 \section{Vanishing loci of modular forms.}

Our aim is to show the smoothness and to determine the image under the building map $\la$ of vanishing sets of modular
forms. Obviously, it suffices to do so in the fundamental domain $\bF \subset \OM^r$ for $\Ga$. Therefore we define
for a modular form $f$ the analytic space
 $$V(f) := \{\bo \in \bF~|~ f(\bo) = 0\}.\leqno{(4.1)}$$
We start with $f = g_i$, where $1 \leq i < r$ ($g_r = \De$ never vanishes on $\OM^r$). It has been shown in \cite{10},
Corollary 3.6 that $V(g_i)$ is contained in $\bF_{r-i} = \{\bo \in \bF~|~ |\om_{r-i}| = |\om_{r-i+1}|\}$, that is,
$\la(V(g_i)) \subset \MW_{r-i}(\Q)$. Besides that inclusion, the size $|g_i(\bo)|$ for $\bo \in \bF \setminus \bF_{r-i}$ as well
as the spectral norm $\|g_i\|_{\bx}$ for $\bx \in \MW_{r-1}(\Q)$ has been determined in \cite{11}, Corollary 4.16. Here we show:
 \vspace{0.3cm}

{\bf 4.2 Theorem.} {\it (i) For $1 \leq i < r$, $\la(V(g_i)) = \MW_{r-i}(\Q)$. More generally, let $S$ be any non-empty subset of 
$\{1,2,\ldots,r-1\}$. Then
 $$\la(\bigcap_{i\in S} V(g_i)) = \bigcap_{i\in S} \MW_{r-i}(\Q).$$
(ii) The $V(g_i)$ ($i\in S$) intersect transversally, and the analytic space $\underset{i\in S}{\bigcap} V(g_i)$ is smooth 
of dimension $r-1-\#(S)$.}
 \vspace{0.3cm}

In what follows, we will at various places make use of the
 \vspace{0.3cm}

{\bf 4.3 Observation.} Knowledge of the following data on $\bo = (\om_1,\ldots,\om_{r-1},1) \in \bF$ is equivalent:
 \begin{itemize}
 \item[(a)] $|\om_1|, \ldots,|\om_{r-1}|$;
 \item[(b)] $|\mu_1|,\ldots,|\mu_r|$, where $\mu_i = \mu_i(\bo)$;
 \item[(c)] $\bx = \la(\bo)$;
 \item[(d)] the Newton polygon $NP_{\bo}$ of $T^{-1}\phi_T^{\bo}(X)$.
 \end{itemize}

 \begin{proof}
First note that $\{\mu_r, \ldots,\mu_1\}$ is an $\F$-SMB of $V = {_T\phi^{\bo}}$ (\cite{10}, 3.4), so (b) determines the
spectrum of $V$. Now (a) $\Rightarrow$ (b) is \cite{11}, 4.2, (b) $\Leftrightarrow$ (d) is (1.10), (a) $\Leftrightarrow$ (c)
is trivial, as $x_i = \log\,\om_i$. A closer look to the equations \cite{11}, 4.2 shows that they may be solved for the $|\om_i|$ if
the $|\mu_i|$ are given, thus (b) $\Rightarrow$ (a). 
 \end{proof}

 \begin{proof} of 4.2:
(i) The inclusion of the left hand side into the right hand side is trivial. For the reverse inclusion, let $\bx \in \underset{i\in S}
{\bigcap} \MW_{r-i}(\Q)$ and $\bo \in \bF_{\bx}$. Then $|\om_{r-i}| =|\om_{r-i+1}|$ for $i\in S$, thus also 
$|\mu_{r-i}| = |\mu_{r-i+1}|$ (\cite{10}, 3.4). The remark in (1.11) along with (3.2) shows that there is an $\bo'$ with 
$NP_{\bo} = NP_{\bo'}$ such that $g_i(\bo') = 0$ for all $i \in S$. (Note the reverse order in the $\F$-SMB $\{\mu_r,\ldots,\mu_1\}$!)
So $\bx = \la(\bo) = \la(\bo')$ lies in the left hand side.
 \medskip

(ii) This follows from the non-vanishing of $\underset{1 \leq i,j \leq r-1}{\det} (\frac{\partial g_i}{\partial \om_j})$,
cf. (3.15).
 \end{proof}

{\bf 4.4 Example.} Let $S$ be the full set $\{1,2,\ldots,r-1\}$. For $\bo \in \bF$, the following are equivalent:
 \begin{itemize}
 \item[(a)] $g_i(\bo) =0 \quad \forall i \in S$;
 \item[(b)] $\alpha_i(\bo) = 0 \quad \forall i \in S$;
 \item[(c)] $E_{q^i-1}(\bo) = 0 \quad \forall i \in S$;
 \item[(d)] $\bo \in \OM^r(\F^{(r)})$.
 \end{itemize}
Here (a) $\Leftrightarrow$ (b) $\Leftrightarrow$ (c) follows from (3.4), and its equivalence with (d) is \cite{10},
Proposition 2.9. Correspondingly, $\la(\underset{i \in S}{\bigcap} V(g_i)) = \la(\OM^r(\F^{(r)})) = \{\boldsymbol 0\}$.
 \vspace{0.3cm}

Next, we deal with the Eisenstein series $E_k(\bo)$, where we always suppose that $0 < k \equiv 0 \,(\bmod\,q-1)$, as 
$E_k$ vanishes identically if $k \not\equiv 0\,(\bmod\,q-1)$.
 \vspace{0.3cm}

{\bf 4.5 Theorem.} {\it (i) $\la(V(E_k)) = \MW_{r-1}(\Q)$;
 \medskip

(ii) $V(E_{q^j-1})$ is smooth for $j > 0$;
 \medskip

(iii) the same statement as (4.2)(ii), with $g_i$ replaced by $E_{q^i-1}$.}
 
 \begin{proof}
(iii) follows from the non-vanishing of functional determinants as in (4.2).
 \medskip

(i) We have
 $$E_k(\bo) = \underset{\ba\in A^r}{\sum{}^{'}} (a_1\om_1+ \cdots + a_r\om_r)^{-k}.\leqno{(4.5.1)}$$
Let $i$ be minimal with $|\om_i| = 1$. Then $i < r \Leftrightarrow \bo \in \bF_{r-1}$. The terms of largest absolute value in
(4.5.1) are the $(a_i\om_i+ \cdots + a_r\om_r)^{-k}$ with $a_j \in \F$, all of absolute value 1. If $i=r$, they sum up to
$\underset{a\in \F^*}{\sum} (a\om_r)^{-k} = -1$, so $E_k(\bo) \not= 0$ and $\la(V(E_k)) \subset \MW_{r-1}$. Conversely, let
$\bx \in \MW_{r-1}$, so $i < r$ for $\bo \in \bF_{\bx}$. There do exist $\om_i,\ldots,\om_r \in \overline{\F}$ linearly independent over
$\F$ such that $\underset{a_i,\ldots,a_r\in \F}{\sum{}^{'}} (a_i\om_i+ \cdots+ a_r\om_r)^{-k} = 0$, since the vanishing set of
this rational function on $\overline{\F}^{r-i+1}$ cannot be contained in the union of the $\F$-rational hyperplanes. 
Therefore, the canonical reduction of $E_k(\bo)$ as a function on the affinoid $\bF_{\bx}$ has zeroes and $E_k(\bo)$ has non-constant
absolute value on $\bF_{\bx}$. By Theorem 2.4, $E_k$ presents a zero on $\bF_{\bx}$, which shows that $\bx \in \la(V(E_k))$.
 \medskip

(ii) For $j < r$, the smoothness of $V(E_{q^j-1})$ is covered by (iii). In the general case, let $\bo \in \bF_{r-1}$ be a
zero of $E_{q^j-1}$ and $i < r$ as in (i), $D = \frac{\partial}{\partial \om_i}$. Then
 $$\begin{array}{lll}
D E_{q^j-1}(\bo) &=& {\displaystyle \underset{\ba\in A^r}{\sum}{}^{'} a_i(a_1\om_1+\ldots+ a_r\om_r)^{-q^j}}
 \vspace{0.2cm}\\
 &\equiv& {\displaystyle \underset{a_i,\ldots,a_r\in \F}{\sum{}^{'}}a_i(a_i\om_i,\ldots,a_r\om_r)^{-q^j}}
 \end{array}$$
(where ``$\equiv$'' means up to terms of value strictly less than 1)
 $$ \equiv -\sum_{a_{i+1},\ldots a_r\in \F} (\om_i+a_{i+1}\om_{i+1}+ \ldots+ a_r \om_r)^{-q^j}.$$
Let $V$ be the $\F$-vector space spanned by $\om_{i+1},\ldots,\om_r$, with $e$-function $e_V$. Then
 $$e_V(z)^{-1} = \sum_{v\in V} \frac{1}{z+v},$$
and the above is $-e_V(\om_i)^{-q^j}$, of absolute value 1. Hence $D E_{q^j-1}(\bo) \not= 0$, which shows
the smoothness in $\bo$.
 \end{proof}

{\bf Remark.} In the proof of (i) we have seen that $E_k(\bo) \equiv -1$ if $\bo \in \bF \setminus \bF_{r-1}$,
i.e., if $\la(\bo) \not\in \MW_{r-1}(\Q)$, while $\|E_k\|_{\bx} = 1$ for $\bx \in \MW_{r-1}(\Q)$. Concerning
smoothness and intersection properties of their vanishing sets, in contrast with the ``special'' Eisenstein
series of weight $q^i-1$ for some $j$, the ``non-special'' ones behave rather erratic and unpredictably. 
This holds already in the case $r=2$, see e.g. \cite{8}, Remark 6.7.
 \vspace{0.3cm}

In order to treat the para-Eisenstein series $\alpha_k$, we make the following definition.
 \vspace{0.3cm}

{\bf 4.6 Definition.} {\it For $k \in \N$ let
 $$\bF(k) := \{\bo \in \bF~|~ \La_{\bo} \mbox{ is $k$-inseparable}\}.$$
As this depends only on $\la(\bo)$, we also put
 $$\begin{array}{ll}
 \MW(k)(\Q) := \la(\bF(k)) = \{\bx \in \MW(\Q)~|& \La_{\bo} \mbox{ is $k$-inseparable for}\\
 & \mbox{one (thus all) } \bo \in \la^{-1}(\bx)\}.
 \end{array}$$
It is the set of $\Q$-points of a full simplicial subcomplex $\MW(k)$ of $\MW$.}
 \vspace{0.3cm}

{\bf 4.7 Examples.} Let $\{\la_1,\la_2,\ldots\}$ be the $\F$-SMB of $\La_{\bo}$ constructed in (1.12).
Hence $\bo \in \bF(k)$ if and only if $|\la_k| = |\la_{k+1}|$.
 \medskip

(i) $\MW(1) = \MW_{r-1}$, for $|\la_1| = |\la_2| \Leftrightarrow |\om_{r-1}| = |\om_r| = 1$.\\
Assume now that $r=3$. With some labor, we find the following descriptions of $\MW(2)$, $\MW(3)$ ,
$\MW(4)$ through their vertices. Note that the vertices of $\MW$ are the $[L_{\bk}]$, where $\bk$ is an
$\N_0$-combination of the $\bk_i = (1,1,\ldots,1,0,\ldots,0)$ ($i$ ones, $1 \leq i < r$).
 \medskip

(ii) vertices of $\MW(2)$: $\bk = \boldsymbol 0$ and $\bk = \bk_2+n\bk_1$ ($n \in \N_0$) (see (1.12)).
 \medskip

(iii) vertices of $\MW(3):\: \bk_2$, $2\bk_2+n\bk_1$ ($n \in \N_0$), $n\bk_1$ ($n \in \N$). This corresponds 
to the fact
that $|\la_3| = |\la_4|$ if either of the conditions is satisfied:
 \begin{itemize}
 \item $|\om_2| = |T^2\om_3| = q^2$
 \item $q \leq |\om_1| = |\om_2| \leq q^2$
 \item $|\om_1| = q$
 \item $|\om_2| = |\om_3| = 1$ and $|\om_1| \geq q$.
 \end{itemize}

(iv) vertices of $\MW(4)$: ${\bf 0},\bk_1,\bk_1+\bk_2+n\bk_1$ ($n \in \N_0$), $2\bk_2$, $3\bk_2+n\bk_1$ ($n \in \N_0$).
This corresponds to: $|\la_4| = |\la_5|$ if and only if one of the following holds:
 \begin{itemize}
 \item $|\om_2| = |\om_3| = 1$ and $|\om_1| \leq q$
 \item $|T\om_2| = |\om_1| \leq q^2$
 \item $|\om_2| = |T\om_3| = q$ and $|\om_1| \geq q^2$
 \item $|\om_1| = |T^2\om_3| = q^2$ and $q \leq |\om_2| \leq q^2$
 \item $q^2 \leq |\om_1| = |\om_2| \leq q^3$.
 \medskip\end{itemize}

In the examples, $\MW(k)$ is a full subcomplex of $\MW$ which is everywhere of dimension $r-2$ (that is,
each simplex belongs to a simplex of maximal dimension $r-2$), connected and contractible. These properties should hold
in full generality. Certainly, the $\MW(k)$ deserve more investigation!
 \vspace{0.3cm}

{\bf 4.8 Theorem.} {\it
(i) $\la(V(\alpha_k)) = \MW(k)(\Q)$;
 \medskip

(ii) the same statement as $(4.2)${\rm (ii)} with $g_i$ replaced by $\alpha_i$.}
 
 \begin{proof}
Once again, (ii) follows from the non-vanishing of the functional determinant.
 \medskip

(i)(a) The fact that $V(\alpha_k)$ is contained in $\bF(k)$ is a consequence of (1.11) and the definition
of $\bF(k)$. Suppose that $\bx \in \MW(k)(\Q)$. We will show that $|\alpha_k|$ is not constant on
$\bF_{\bx}$, which forces the existence of a zero of $\alpha_k$ in $\bF_{\bx}$ and gives the wanted 
equality.
 \medskip

(b) From $e_{\bo}(X) = \underset{k\geq 0}{\sum} \alpha_k(\bo) X^{q^k} = X\underset{\la\in\La_{\bo}}
{\prod{}^{'}} (1-X/\la)$ we see that
 $$\alpha_k = \alpha_k(\bo) = s_{q^k-1} \{\la^{-1}~|~0 \not= \la \in \La_{\bo}\}$$
is the $(q^k-1)$-th elementary symmetric function in the $\la^{-1}$:
 $$ \alpha_k = \sum_{S} P(S),$$
where $S$ runs through the family of $(q^k-1)$-subsets of $\La_{\bo} \setminus \{0\}$ and 
 $$P(S) = (\prod_{\la \in S} \la)^{-1}.$$
(c) Let $m+1$ (resp. $n$) be the least (resp. largest) subscript $i$ such that $|\la_i| = |\la_k| = |\la_{k+1}|$,
where $\{\la_1,\la_2,\ldots\}$ is an $\F$-SMB of $\La_{\bo}$ formed out of the $\{T^j\om_i\}$ (see (1.12)). Then
$m < k < n$, the $\om_i$ appearing in $\la_{m+1},\ldots,\la_n$ are all different, and therefore $n-m \leq r$.
 \medskip

(d) Some $P(S)$ has largest absolute value if $S$ contains all the $q^m-1$ elements of $V' \setminus\{0\}$ and
$q^k-q^m$ many elements of $V \setminus V'$, where $V=\underset{1\leq i\leq n}{\sum} \F\la_i$, 
$V' = \underset{1 \leq i \leq m}{\sum} \F\la_i$. The contribution of such $S$ to $\alpha_k$ is
 $$P := (\underset{\la \in V'}{\prod}\la^{-1} ) \underset{S' \subset V\setminus V' \atop\#(S') =q^k-q^m}{\sum}
  P(S'),\: P(S') = (\underset{\la \in S'}{\prod} \la)^{-1}.$$
All the $P(S)$ of such $S$ have the same absolute value $|P(S)|=: c$, which depends only on $|\la_1|,\ldots,|\la_n|$ 
and therefore only on $\bx$. We write $x \equiv y$ if $|x-y| < c$. Then
 $$\alpha_k(\bo) \equiv P \equiv \alpha_k(V),$$
all of which are homogeneous functions of $\bo \in \bF_{\bx}$ of weight $q^k-1$. (That is, each of these
functions $f$ on $\bF_{\bx}$ may be regarded as a homogeneous function $f^*$ on the cone 
$\bF_{\bx}^* = \{(\om_1,\ldots,\om_r) \in C_{\infty}^r~|~ (\om_1 :\cdots: \om_r) \in \bF_{\bx}\}$, and
$f^*(t\bo) = t^{1-q^k}f^*(\bo)$.)
 \medskip

(e) Therefore we must show that $|\alpha_k(V)|$ is not constant on $\bF_{\bx}$, which, as 
$\alpha_n(V) = (\underset{\la \in V}{\prod'}\la)^{-1}$ has constant absolute value, is equivalent 
with the fact that the $k$-th coefficient $\ga_k(V)$ of 
 $$f_V(X) := \alpha_n(V)^{-1} e_V(X) = \prod_{\la \in V}(X-\la) = \sum_{0 \leq k \leq n} \ga_k(V) X^{q^k}$$
has non-constant absolute value. This will be shown by varying the leading coefficients in $\overline{\F}$
of the $n-m$ different $\om_i$ that appear in $\la_{m+1},\ldots,\la_n$ (see (c)) without changing 
$\bx = \la(\bo)$. Dropping the requirement $\om_r =1$, we rescale  the projective coordinates of
$\bo$ such that $|\la_j| < 1$ for $1 \leq j \leq m$ and $|\la_j| = 1$ for $m+1 \leq j \leq n$. Then the
monic polynomial $f_V$ gets coefficients in $O_{C_{\infty}}$, and its reduction $\overline{f}_V \in
\overline{\F}[X]$ satisfies
 $$(\overline{f}_V) = (f_{\overline{V}})^{q^m},$$
where $\overline{V} \subset \overline{\F}$, spanned by $\overline{\la}_{m+1},\ldots,\overline{\la}_n$, is the
reduction of $V$ and $f_{\overline{V}}(X)=  \underset{\overline{\la} \in \overline{V}}{\prod} (X-\overline{\la})$.
 \medskip

(f) Let $\ell := n-m = \dim\,\overline{V}$, which is larger or equal to 2. For an $\F$-subspace $U$ of 
$\overline{\F}$ of dimension $\ell$, write
 $$f_U(X) = \prod_{u\in U} (X-u) = \sum_{0 \leq j \leq \ell} \ga_j(U)X^{q^j}.$$
The $\ga_j = \ga_j(U)$ are homogeneous of weight $(q^j-q^{\ell})$ (that is, $\ga_j(tU) = t^{q^{\ell}-q^j}\ga_j(U),\,
t \in \overline{\F}^*$) and may be considered as forms on $\OM^{\ell}(\overline{\F}) = \{\bnu = (\nu_1: \cdots: \nu_{\ell})
\in \BP^{\ell-1}(\overline{\F})~|~ \nu_1,\ldots,\nu_{\ell} \: \F{\rm -l.i.}\}$. Further, for $0 < j < \ell$, $\ga_j(\bnu)$ 
vanishes somewhere on $\OM^{\ell}(\overline{\F})$: take e.g. the $\nu_j$ as an $\F$-basis of $U := \F^{(\ell)}$; 
then $f_U(X) = X^{q^{\ell}}-X$.
 \medskip

(g) Putting $U := \overline{V}$, of dimension $\ell = m-n$, $j := k-m$, the fact that $\ga_j(\overline{V})$
takes both zero and non-zero values on $\OM^{\ell}(\overline{\F})$ implies that $|\alpha_k(V)|$ is
non-constant on $\bF_{\bx}$, provided we can find these zero/non-zero values inside the image of
$\bF_{\bx}$ in $\OM^{\ell}(\overline{\F})$.
 \medskip

(h) Let $I$ be the set of those indices $i$ with $1 \leq i \leq r$ and $\log\,\om_i \in \Z$. Due to our
normalization $|\la_{m+1}| = \cdots = |\la_n| = 1$, $I' :=
\{m+1,\ldots,k,k+1,\ldots,n\}\subset I$. For
$i \in I$ let $\theta_i \in \overline{\F}^*$ be the leading coefficient of $\om_i$, that is 
 $$\om_i = \theta_iT^{\log\,\om_i} (1+n_i) \quad \mbox{ with } |n_i| < 1.$$
Given $\bo = (\om_1:\cdots: \om_r) \in \bF_{\bx}$, with the scaling of projective coordinates as in (e),
let $\bo'= \bo'_1: \cdots:\bo'_r)$ be defined by
 $$\begin{array}{llll}
 \om'_i &=& \om_i, & \mbox{ if } i \not\in I\\
 &=& \frac{\nu_i}{\theta_i} \om_i, & \mbox{ if } i \in I
 \end{array}$$ 
with some $\nu_i \in \overline{\F}^{\ast}$. Then $\bo'$ belongs to $\bF_{\bx}$ at least if the $\nu_i$ ($i\in I$) are $\F$-linearly 
independent. Choosing first the $\nu_i$ ($i\in I'$) such that
$(\nu_{m+1}:\cdots : \nu_n) \in \OM^{\ell}(\overline{\F})$ is a zero of $\ga_{k-m}$ and then the remaining 
$\nu_i$ ($i \in I\setminus I'$) such that $\{\nu_i~|~i\in I\}$ is $\F$-l.i., we find some $\bo'$ where
$|\alpha_k(V)|$ is strictly less than the spectral norm of $\bo \longmapsto \alpha_k(V)$
on $\bF_{\bx}$. \\
This finishes the proof.
 \end{proof}

Let now $a \in A$ have positive degree, and consider the forms $_a\ell_k$. Similar arguments as in the
proof of Theorem 4.8, based on Proposition 1.11, allow to describe the set $\la(V(_a\ell_k))$, that is,
the analogue of 4.8(i).  However, as both the statement and its proof are substantially more complex than
its counterpart for $\alpha_k$ (see the relatively simple case of $r=2$ treated in \cite{9}, Section 5),
we restrict to stating the analogue of (4.2)(ii), (4.5)(iii) and (4.8)(ii), with identical proof.
 \vspace{0.3cm}

{\bf 4.9 Theorem.} {\it
Let $S$ be a non-empty subset of $\{1,2,\ldots,r-1\}$. The $V(_a\ell_i)$ ($i\in S$) intersect transversally,
and the analytic space $\underset{i\in S}{\bigcap} V(_a\ell_i)$ is smooth of dimension $r-1 = \#(S)$.}
 \hspace{5.1cm} $\Box$
 \vspace{0.3cm}

We conclude with the relationship between $\alpha_k$ and the $_a\ell_k$.
 \vspace{0.3cm}

(4.3) Let $\rho$ be the Carlitz module, i.e., the rank-one Drinfeld module defined by $\rho_T(X) = TX+X^q$, and 
let $_ac_k$ be its coefficients:
 $$\rho_a(X) = \sum_{0 \leq k \leq d} {_ac_k} X^{q^k} \quad (a\in A \mbox{ of degree } d,\: {_ac_0} = a).\leqno{(4.10.1)}$$
We have
 $$\log\,{_ac_k} = (d-k)q^k,\: 0 \leq k \leq d,\: {_ac_k} = 0 \mbox{ for } k > d.\leqno{(4.10.2)}$$
There exists $\overline{\pi} \in C_{\infty}$, well-defined up to a $(q-1)$-th root of unity, such that $\rho = \phi^{(L)}$ is
the Drinfeld module corresponding to the rank-one lattice $L = \overline{\pi} A$, with 
 $$\log\,\overline{\pi} = q/(q-1).\leqno{(4.10.3)}$$
The exponential function of $L$ is
 $$\begin{array}{lll}
 e_L &= &\underset{k\geq 0}{\sum} D_k^{-1} \tau^k,\: D_k:= [k][k-1]^q \cdots [1]^{q^{k-1}},\\
\lbrack k\rbrack  &:=& \lbrack T,k \rbrack = (T^{q^k}-T).
 \end{array}\leqno{(4.10.4)}$$
Hence
 $$\log\,\alpha_k(L) = -kq^k.\leqno{(4.10.5)}$$
Replacing $L$ with $A$ yields the Drinfeld module $\phi^{(A)}$ with 
 $$\phi_T^{(A)} = TX+ \overline{\pi}^{q-1}X^q$$
and
$$\begin{array}{cll}
 \alpha_k(A) &=& \overline{\pi}^{(q^k-1)} \alpha_k(L),\\
 \log\,\alpha_k(A) &=& q(q^k-1)/(q-1)-kq^k.
 \end{array}\leqno{(4.10.6)}$$
All of this is easily verified and may be found at different places, e.g. \cite{5}, \cite{13}, \cite{16}. We also
need the following lemma, whose proof is an exercise in manipulating the preceding formulas.
 \vspace{0.3cm}

{\bf 4.11 Lemma.} {\it
If $d = \deg\,a$ tends to infinity then $D_k{_ac_k}/\lbrack a,k\rbrack$ tends to $1$.}
 \vspace{0.3cm}

Now we normalize the modular forms $_a\ell_k$ and $\alpha_k$ by dividing through the corresponding quantities of
the Drinfeld module $\phi^{(A)}$.
 $${_a\tilde{\ell}_k} := \overline{\pi}^{1-q^k} {_ac_k}^{-1} {_a\ell_k}; \: \tilde{\alpha}_k := \overline{\pi}^{1-q^k}D_k \alpha_k.
 \leqno{(4.12)}$$
This normalization is quite natural, see Remark 4.14.
 \vspace{0.3cm}

{\bf 4.13 Theorem.} {\it
As the degree $d$ of $a$ tends to infinity, $_a\tilde{\ell}_k$ tends to $\tilde{\alpha}_k$, locally
uniformly on $\OM^r$.}
 \medskip

(Here locally uniform convergence means uniform convergence on the parts of an admissible covering of
$\OM^r$.)
 
 \begin{proof}
(i) Let $V$ be the space of isobaric polynomials of weight $q^k-1$ in the $g_1,\ldots,g_r$ (which by (1.8)
is the space of modular forms of weight $q^k-1$ and type 0). As $\dim_{C_{\infty}}(V) < \infty$, all norms on
$V$ agree. Since $_a\ell_k$, $\alpha_k$ and their normalizations belong to $V$, it suffices to show convergence
with respect to one specific norm on $V$.
 \medskip

(ii) We let $\bF^{(k)} := \{\bo \in \bF~|~\log \om_{r-1} \geq k\}$. It is an open admissible subspace, and we will
use the norm $\|f\| := \underset{\bo \in \bF^{(k)}}{\sup} f(\bo)$ (well-defined in view of (1.8)(c)).
 \medskip
For $\bo \in \bF$ let $\{\la_1(\bo),\la_2(\bo), \ldots \}$ be the $\F$-SMB of $\La_{\bo}$ as in (1.12). Now
if $\bo \in \bF^{(k)}$ then for $1 \leq i \leq k$, $\la_i(\bo) = T^{i-1}\om_r = T^{i-1}$. Hence the
$\F$-span of these agrees with $A_{k-1} := \{a \in A~|~ \deg\,a \leq k-1\}$. As in the proof of 4.8 we
find for $1 \leq i \leq k$
 $$\alpha_i(\bo) = (\underset{c\in A_{i-1}}{\prod{}^{'}}c)^{-1} + \mbox{ smaller terms},$$
that is
 $$|\alpha_i(\bo)| = |\alpha_i(A)| = |\overline{\pi}|^{q^i-1} |D_i|^{-1}.$$
(iii) Suppose that $\bo \in \bF^{(k)}$ and $d = \deg\,a \geq k$. The elements  
 $$\mu_{i,j} := e_{\bo}(T^{j-1}\om_i/a) \mbox{ with } 1 \leq i \leq r,\: 1 \leq j \leq d$$
form an $\F$-basis of $_a\phi^{\bo}$, the $a$-division points of $\phi^{\bo}$. We see from
 $$
 \mu_{i,j}(\bo)| = |\frac{T^{j-1}\om_i}{a}| \underset{\la \in \La_{\bo}\atop |a\la|<|T^{j-1}\om_i|}{\prod{}^{'}}
|1 -\frac{T^{j-1}\om_i}{a\la}|$$
that the $\la_j:= \mu_{r,j}$ with $1 \leq j \leq k$ are the first $k$ elements of an $\F$-SMB of $_a\phi^{\bo}$. Moreover,
the next $\F$-SMB vector $\la_{k+1}$ satisfies $|\la_{k+1}| > |\la_k|$, due to the assumption $\bo \in \bF^{(k)}$. From
 $$\phi_a^{\bo}(X) = \sum_{0 \leq i \leq rd} {_a\ell_i} (\bo) X^{q^i}$$
we find
 $$\begin{array}{lll}
 a^{-1}{_a\ell_i}(\bo) &=& \alpha_i(_a\phi^{\bo}) = s_{q^i-1} \{\mu^{-1}~|~ 0 \not= \mu \in {_a\phi^{\bo}}\}
 \vspace{0.3cm}\\
 &=& (\underset{c\in A_{i-1}}{\prod{}^{'}}(\frac{c}{a}) \cdot u \end{array}$$ 
with some $u \in C_{\infty}$ of absolute value 1. This implies
 $$|{_a\ell_i}(\bo)| = |{_a\ell_i}(A)| = |\overline{\pi}|^{q^i-1} |{_ac_i}|,$$
valid for $i \leq k \leq d$ and $\bo \in \bF^{(k)}$.
 \medskip

(iv) Next, consider the identity (3.4)(ii)
 $$[a,k] \alpha_k(\bo) = \sum_{1\leq i \leq k-1} {_a\ell_i}(\bo) \alpha_{k-i}^{q^i} (\bo) + {_a\ell_k}(\bo).$$
For $d = \deg\,a \geq k$, all the terms have constant absolute value on $\bF^{(k)}$. Plugging in, we see that
$\log([a,k]\alpha_k)$ and $\log({_a\ell_k})$ grow of order $(d-k)q^k+q(q^k-1)/(q-1)$ with $d \lra \infty$, while the log
of the other terms grow of order less or equal to $(d-k+1)q^{k-1}+q(q^k-1)/(q-1)$. Upon normalization 
$f \rightsquigarrow \tilde{f}$, we find that ${_a\tilde{\ell}_k}$ tends to $[a,k]/({_ac_k} D_k)^{-1} \cdot \tilde{\alpha}_k$
uniformly on $\bF^{(k)}$. The result now follows from Lemma 4.11. 
 \end{proof}

{\bf 4.14 Remarks.} (i) In steps (ii) and (iii) of the preceding proof, in fact the stronger statements hold:
 $$\lim\,\alpha_i(\bo) = \alpha_i(A),\: \lim\,{_a\ell_i} (\bo) = {_a\ell_i}(A),$$
where the limits are with respect to $|\om_{r-1}| \lra \infty$. This follows from a closer look to the arguments and
estimates used there.  Hence $\lim\,\tilde{\alpha}_i (\bo) = 1 = \lim\,{_a\tilde{\ell}_i} (\bo)$, which in the case
$r=2$ means $\tilde{\alpha}_i(\infty) = 1 = {_a\tilde{\ell}_i} (\infty)$.
 \medskip

(ii) Theorem 4.13 has been shown in the case $r=2$ in \cite{9}, Theorem 6.16. The present proof isn't but a generalization
of this special case.
 \vspace{0.3cm}

{\bf 4.15 Concluding remarks/questions.}
 \medskip

(i) Let $f$ be one of the functions $E_k$ or $\alpha_i$ on $\OM^r$. For $\bx \in \la(V(f))$, $f$ is given on $\bF_{\bx}$ as a convergent
sum of terms 
 $$f(\bo) = \sum T_{\bsi} + \sum U_{\bsj}$$
with finitely many $T_{\bsi}$ all of constant absolute value $|T_{\bsi}| = \|f\|_{\bx}$ and terms $U_{\bsj}$ of strictly smaller
value.
 \medskip

For $f = E_k$, the $T_{\bsi}$ are the $|a_i\om_i+ \cdots + a_r\bo_r|^{-k}$, $a_i,\ldots,a_r \in \F$, see proof of (4.5); for
$f = \alpha_i$, the $T_{\bsi}$ are certain $P(S)$, see proof of (4.8), part (d). A similar property may be shown for 
$f = {_a\ell_k}$. Hence $|f(\bo)| < \|f\|_{\bx}$ arises from cancellations between the terms $T_{\bsi}$. The proof of
(2.6) combined with Remark 2.9 yields the following generalization of (2.6):
 \medskip

Given any of the modular forms $f$ as above, the map
 $$ \bx \longmapsto \log_q\|f\|_{\bx}$$
is an affine function on $\MB\MT(\Q)$ (which off $\la(V(f))$ agrees with $\log\,f(\bo)$, $\bo \in \bF_{\bx}$).
 \medskip

To which class of modular forms does this property generalize?
 \medskip

(ii) For all our distinguished modular forms $f$ ($f = \alpha_i,{_a\ell_j},E_k$), $\la(V(f))$ is (the set of
$\Q$-valued points of) a simplicial subcomplex of pure dimension $r-1$ of $\bF$. In the case $r=2$ this means that
the zeroes of $f$ lie in the $\la$-preimage of vertices of the Bruhat-Tits tree $\MB\MT$ of ${\rm PGL}(2,K_{\infty})$. 
 \smallskip

How can we characterize modular forms with this property?

\end{document}